\documentclass[11pt,leqno]{amsart}
\usepackage[a4paper,margin=28mm]{geometry}
\usepackage{amsmath,amssymb,amsthm,mathtools}
\usepackage[T1]{fontenc}
\usepackage{lmodern}
\usepackage{microtype}
\usepackage{needspace,booktabs}
\usepackage[colorlinks=true,linkcolor=blue,citecolor=blue,urlcolor=blue]{hyperref}
\hypersetup{pdftitle={Euler Trace Congruences over Local Fields},
 pdfauthor={Wei Tao and Guo Xuejun},
 pdfsubject={Exact Euler trace congruences for prescribed residual similarity classes}}
\allowdisplaybreaks
\numberwithin{equation}{section}
\newtheorem{theorem}{Theorem}[section]
\newtheorem{proposition}[theorem]{Proposition}
\newtheorem{lemma}[theorem]{Lemma}
\newtheorem{corollary}[theorem]{Corollary}
\theoremstyle{definition}

\newtheorem{example}[theorem]{Example}
\theoremstyle{remark}
\newtheorem{remark}[theorem]{Remark}
\newcommand{\OO}{\mathcal O}
\newcommand{\QQ}{\mathbb Q}
\newcommand{\ZZ}{\mathbb Z}
\newcommand{\FF}{\mathbb F}
\newcommand{\barA}{\overline A}
\newcommand{\vp}{v_\pi}
\newcommand{\Sps}{S}
\DeclareMathOperator{\Tr}{Tr}
\DeclareMathOperator{\Mat}{M}
\DeclareMathOperator{\ord}{ord}

\title[Euler trace congruences over local fields]{Euler Trace Congruences over Local Fields}
\author{Wei Tao}
\author{Guo Xuejun}
\date{}
\begin{document}
\begin{abstract}
	Let $K$ be a finite extension of $\QQ_p$ with uniformizer $\pi$
	and residue field $\FF_q$. For each $r\ge1$, we determine the
	largest integer $E_r$ such that
	$\Tr(A^{q^r})\equiv\Tr(A^{q^{r-1}})\pmod{\pi^{E_r}}$
	for every integral matrix $A$ with a prescribed residual similarity class.
	We give formulas for both nilpotent and nonnilpotent residual
	classes and construct integral lifts attaining the resulting
	exponents.
	We also exhibit local fields with the same ramification index
	and residue field for which the optimal exponents differ,
	even when $r$ and the residual similarity class are fixed.
\end{abstract}
\maketitle

\section{Introduction}
For an integer matrix $A$, a prime $p$, and a positive integer $r$,
the $r$th matricial Euler congruence states that
$\Tr(A^{p^r})\equiv\Tr(A^{p^{r-1}})\pmod{p^r}$.
This is the case $k=p^r$ of the Gauss congruences
$\sum_{d\mid k}\mu(k/d)a_d\equiv0\pmod k$
for the sequence $a_k=\Tr(A^k)$, where $\mu$ is the
M\"obius function.

The same congruences hold when $a_k$ counts the fixed points
of the $k$th iterate of a map: the sum
$\sum_{d\mid k}\mu(k/d)a_d$, divided by $k$, counts orbits
of length $k$.
For the trace sequence, the identity
\[
\exp\left(\sum_{k\ge1}\frac{\Tr(A^k)}{k}t^k\right)
=\det(I-tA)^{-1}
\]
shows that the series on the left has integer coefficients.
Byszewski, Graff, and Ward \cite{BGW} survey these connections
between Dold congruences, traces, and periodic points.

Zarelua \cite{Zarelua} studied congruences for traces of powers and their
topological and dynamical interpretations. The historical discussion in Mazur and Petrenko
\cite{MazurPetrenko} traces the matricial Euler congruence
back to J\"anichen and Schur, before its renewed study by
Arnold, Zarelua, and Vinberg.

Deligne \cite{Deligne} established an extension to truncated Witt vectors,
with Frobenius acting on the trace of the preceding power.
Mazur and Petrenko \cite{MazurPetrenko} proved that congruence modulo $p^j$
of the characteristic polynomials of two integer matrices implies
congruence modulo $p^{j+1}$ for their $p$th powers.
Steinlein \cite{Steinlein} gave elementary matrix versions of the Fermat
and Gauss congruences. More recently, Abyzov \cite{Abyzov} extended the Euler and
Gauss congruences to matrices over commutative rings.
His Theorem~4.2 gives a common divisibility bound for all
matrices over a ring with a finite residue field, without
imposing any condition on the residual similarity class.

Over a finite extension of $\QQ_p$, we determine the largest
exponent for which the Euler trace congruence holds for every
integral lift of a prescribed residual similarity class.
For each class, we construct an integral lift attaining the stated exponent.

Let $K/\QQ_p$ be finite, with ring of integers $\OO$, uniformizer $\pi$,
and residue field $\FF_q$, where $q=p^f$. Normalize $\vp(\pi)=1$, set
$e=\vp(p)$ and $\vp(0)=+\infty$, and extend this normalization to
unramified extensions.

Fix positive integers $m$ and $r$, and let
$
D_r(A)=\operatorname{Tr}(A^{q^r})
-\operatorname{Tr}(A^{q^{r-1}})
$
 be the trace difference for $A\in M_m(\mathcal O)$.
For $\overline A\in M_m(\mathbb F_q)$, let
\begin{equation}\label{eq:definition}
	E_r(K,\overline A)
	=
	\min_{\substack{
			A\in M_m(\mathcal O)\\
			A\bmod\pi\sim\overline A
	}}
	v_\pi\bigl(D_r(A)\bigr)
\end{equation}
be the $r$th optimal congruence exponent for the
residual similarity class of $\overline A$.

A conjugating matrix over the residue
field lifts to $\mathrm{GL}_m(\OO)$, so one may equivalently fix a
representative of the residual class. Nilpotence always refers to
$\barA$, not to an integral lift. We write $E_r$ when the field and the
residual class are fixed, and use $I_d$ for the identity matrix of size $d$.
The valuation $v_p$ on integers is normalized by $v_p(p)=1$.

We write $n=f(r-1)$,
$u=p/\pi^e\in\OO^\times$, and
$\eta=\overline u\in\FF_q^\times$. For each positive integer $h$, define
\begin{equation}\label{eq:introL}
	C_n=\min_{0\le i\le n}\{p^i+e(n-i)\},
\ \ 	L_{h,e}(n)=\min_{\lceil p^n/h\rceil\le k\le p^n}
\{k+e(n-v_p(k))\},
\end{equation}
 and $L_{0,e}(n)=+\infty$.
For each positive integer $j$, let $J_j(0)$ denote the
nilpotent Jordan block of size $j$.
For a nilpotent residual class with Jordan partition
$\lambda=(\lambda_1,\ldots,\lambda_t)$, where
$\lambda_1\ge\cdots\ge\lambda_t\ge1$, write
$J_\lambda=\bigoplus_{i=1}^t J_{\lambda_i}(0)$.

\begin{theorem}[Nilpotent residual classes]\label{thm:complete}
For every nilpotent residual class $J_\lambda$ and every $r\ge1$,
we have $E_r=L_{\lambda_1,e}(n)$ except in the three cases
\textup{(I)--(III)} specified in
Section~\ref{sec:exact}.
In those cases the correction is
\begin{equation}\label{eq:answer}
 E_r-L_{\lambda_1,e}(n)=
 \begin{cases}
 1,&\textup{(I)},\ e=1,\\
 \min\{2e,\vp(u^2+1+\pi^{3e/2}u^4)\},&\textup{(I)},\ e\ge2,\\
 \min\{e/3,\vp(u^4+1)\},&\textup{(II)},\\
 \min\{2e/(p-1),\vp(u^2+1)\},&\textup{(III)}.
 \end{cases}
\end{equation}
Every stated value is attained by an integral lift of $J_\lambda$.
\end{theorem}

The largest Jordan block sizes in cases \textup{(I)--(III)}
are $3\cdot2^s$, $5\cdot2^s$, and $(p+1)p^s$,
respectively, for some integer $s\ge0$. The conditions on $K$, on $n=f(r-1)$,
and on the remaining Jordan blocks are listed together at the beginning
of Section~\ref{sec:exact}. In particular, the exceptional correction
can disappear when further Jordan blocks are present.

\begin{corollary}\label{cor:powerblock}
	Suppose that $\barA$ is nilpotent and that its largest Jordan
	block has size $p^s$ for some integer $s\ge0$.
	Then, for every $r\ge1$,
	\begin{equation}
		E_r=\min_{0\le j\le\min(s,n)}\{ej+p^{n-j}\}.
	\end{equation}
\end{corollary}

For a general residual matrix $\barA$, let $\barA_0$ be its
nilpotent primary part and let $\barA_\times$ be its
complementary invertible part, so that
$\barA\sim\barA_0\oplus\barA_\times$.
Let $h_0$ be the largest Jordan block length of $\barA_0$,
with $h_0=0$ if $\barA_0$ is absent.

\begin{theorem}[Nonnilpotent residual classes]\label{thm:nonmain}
	Suppose that $\barA$ is nonnilpotent. Unless there exist
	nonnegative integers $a,b$ such that $f=1$,
	$e=(p-1)p^a$, $n=a+b+1$, and
	$\barA^{p^b-1}
	\bigl(\barA^{p-1}+\eta^{-1}I_m\bigr)^{p^b}=0$,
	we have $E_r=C_n$.
\end{theorem}

Under the exceptional condition in Theorem~\ref{thm:nonmain},
define $\rho_K$ as follows.
For $p>2$, put $c=-\eta^{-1}\in\FF_p^\times$ and let
$\ell=\ord(c)$ be its multiplicative order.
Let $L/K$ be the unramified extension of degree $\ell$.
Choose $\alpha\in\FF_{p^\ell}$ satisfying $\alpha^{p-1}=c$,
and let $\zeta\in\OO_L^\times$ be its Teichm\"uller lift.
Choose $\beta\in\OO_L$ whose residue generates a normal basis
of $\FF_{p^\ell}/\FF_p$; when $\ell=1$, one may take $\beta=1$.
Let
\begin{equation}\label{eq:rho}
	\rho_K=
	\begin{cases}
		\min\!\left\{
		p^a,\,
		\vp\!\left(
		\Tr_{L/K}\bigl(\zeta((1+\pi\beta)^{p^{a+1}}-1)\bigr)
		\right)-p^{a+1}
		\right\},
		& p>2,\\[4pt]
		1,
		& p=2,\ a=0,\\[4pt]
		\min\!\left\{
		3\cdot2^{a-1},\,
		\vp\bigl((1+\pi)^{2^{a+1}}-1\bigr)-2^{a+1}
		\right\},
		& p=2,\ a\ge1.
	\end{cases}
\end{equation}

\begin{theorem}\label{thm:nonexact}
	Suppose that $\barA$ is nonnilpotent and satisfies
	the exceptional condition in Theorem~\ref{thm:nonmain}.
	Then
	\begin{equation}\label{eq:nonanswer}
		E_r=\min\bigl\{C_n+\rho_K,\ L_{h_0,e}(n)\bigr\}>C_n.
	\end{equation}
	The quantity $\rho_K$ is a positive integer independent
	of the choices of $\alpha$, $\beta$, and $\pi$.
	In particular, if $\barA$ is invertible, then
	$E_r=C_n+\rho_K$.
	Every stated value is attained by an integral lift of $\barA$.
\end{theorem}

Section~\ref{sec:prelim} gives preliminary results on integral
lifts, cyclic permutations, and characteristic polynomials.
Sections~\ref{sec:exact} and \ref{sec:nonnil} prove the formulas
for $E_r$ in the nilpotent and nonnilpotent cases, respectively,
and construct lifts attaining these values.
Section~\ref{sec:arithmetic} gives bounds on $\rho_K$ and
examples of local fields with the same $(p,e,f)$ but different
values of $E_r$ for the same residual similarity class.

\section{Preliminary results}\label{sec:prelim}
Lifts with vanishing trace differences will allow us to isolate a
single primary component.

\begin{lemma}\label{lem:periodic}
Every residual similarity class has an integral lift $B$ satisfying
$D_r(B)=0$ for all $r\ge1$.
\end{lemma}
\begin{proof}
For a monic irreducible factor $g\ne T$ of degree $\ell$, take the
unramified extension $L/K$ of degree $\ell$ and a Teichm\"uller lift
$\zeta$ of a root of $g$. The generalized Jordan form in
\cite[Theorem 2.2]{Rybakov} realizes a residual primary block of length
$h$ by restricting scalars from $k_L=\OO_L/\pi\OO_L$ to $\FF_q$ in
$\bar\zeta I_h+J_h(0)$. Lift it to $\zeta I_h+J_h(0)$ over $\OO_L$.
For every positive integer $a$, the binomial expansion and the
vanishing of the traces of positive powers of $J_h(0)$ give
\[
 \Tr((\zeta I_h+J_h(0))^a)=h\zeta^a,\qquad
 \Tr(B^a)=h\Tr_{L/K}(\zeta^a).
\]
Here $B$ denotes the restriction of scalars, and the first trace is
the matrix trace over $L$. If $\sigma$ generates $\operatorname{Gal}(L/K)$
with $\sigma(\zeta)=\zeta^q$, then trace invariance gives
\[
 \Tr_{L/K}(\zeta^{q^r})
 =\Tr_{L/K}(\sigma(\zeta^{q^{r-1}}))
 =\Tr_{L/K}(\zeta^{q^{r-1}}).
\]
Thus this block contributes zero to every $D_r$. For $g=T$, use
$J_h(0)$ itself; every positive power has trace zero. Direct sums preserve
the required residual type and add the trace differences, so the sum
of these lifts has $D_r(B)=0$ for every $r$.
\end{proof}

These lifts allow us to compute $E_r$ separately on the primary parts.

\begin{proposition}\label{prop:primary}
Let $\barA_g$ be the primary parts of $\barA$, indexed by the distinct
monic irreducible factors of its characteristic polynomial. Every
integral lift splits into invariant integral summands lifting these
parts, and
\begin{equation}\label{eq:primarymin}
 E_r(K,\barA)=\min_g E_r(K,\barA_g).
\end{equation}
\end{proposition}
\begin{proof}
Write the residual characteristic polynomial as $\prod_g g^{d_g}$.
Its factors are pairwise relatively prime. Hensel factorization gives
$\det(TI_m-A)=\prod_g F_g(T)$ with $\overline F_g=g^{d_g}$.
The polynomials $F_g$ are pairwise relatively prime over $\OO[T]$:
their resultants reduce to nonzero elements of the residue field and
are units. The Chinese remainder theorem therefore supplies polynomial
idempotents $e_g(T)$ modulo $\prod_g F_g$. By Cayley--Hamilton, the
operators $e_g(A)$ are orthogonal idempotents summing to $I_m$.
Hence
\[
 \OO^m=\bigoplus_g e_g(A)\OO^m.
\]
Each summand is free over the discrete valuation ring $\OO$, is stable
under $A$, and reduces to the prescribed primary subspace.
The trace differences add, so every lift has valuation at least the
right side of \eqref{eq:primarymin}.
Conversely, choose a primary component on which the right side of
\eqref{eq:primarymin} is attained. On this component take a lift
attaining its minimum, and on every other component use
Lemma~\ref{lem:periodic}.
The sum of their trace differences is exactly that of the chosen
component. A nonempty set of nonnegative integral valuations has a
minimum; if all values are infinite, any lift suffices.
This proves the reverse inequality.
\end{proof}

The decomposition above is over distinct irreducible factors.
Different Jordan blocks belonging to the same factor need not split
in an integral lift. To estimate trace differences, we use the
following consequence of cyclic invariance of the trace.

\begin{lemma}\label{lem:words}
For $s\ge0$ and $1\le t\le p^s$, let $H_{s,t}(B,X)$ be the coefficient of $z^t$ in
$\Tr((B+zX)^{p^s})-\Tr(B^{p^s})$. As a polynomial in the entries
of $B$ and $X$, $H_{s,t}$ has integer coefficients divisible by
$p^{\max(s-v_p(t),0)}$. Moreover, for $0\le a\le s$,
\begin{equation}\label{eq:word-leading}
 p^{-(s-a)}H_{s,p^a}(B,X)\equiv
 \Tr\bigl(\overline B^{\,p^{s-a}-1}\overline X\bigr)^{p^a}
 \pmod\pi.
\end{equation}
\end{lemma}
\begin{proof}
We use the cyclic-permutation argument of \cite{Steinlein}, recording
the number of occurrences of $X$. Put $P=p^s$. A word is a formal
ordered product $Z_1\cdots Z_P$ with $Z_i\in\{B,X\}$.
The coefficient $H_{s,t}$ is the sum of traces of the words containing
exactly $t$ occurrences of $X$, for $1\le t\le P$. For matrices over a commutative ring,
$\Tr(VW)=\sum_{i,j}V_{ij}W_{ji}=\Tr(WV)$.
Taking $V=Z_1$ and $W=Z_2\cdots Z_P$ shows that cyclic rotation
preserves the trace. If the shortest period of a word has length
$p^{s-b}$, then it is repeated $p^b$ times and its cyclic orbit has
size $p^{s-b}$. Consequently
\[
 p^b\mid t,\qquad b\le v_p(t),\qquad
 \sum_{W'\in\operatorname{Orb}(W)}\Tr(W')=p^{s-b}\Tr(W).
\]
Each orbit sum is therefore divisible, as an integer polynomial, by
$p^{s-v_p(t)}$. This proves the first assertion, including $t=P$.

Now take $t=p^a$. After division by $p^{s-a}$, every orbit with
$b<a$ still has a factor $p$. For $b=a$, the shortest period contains
one $X$ and $p^{s-a}-1$ copies of $B$. All such periods are rotations
of $B^{p^{s-a}-1}X$, and they give one orbit of words of length $P$.
Its contribution is
\[
 p^{s-a}\Tr((B^{p^{s-a}-1}X)^{p^a}).
\]
For a matrix $M$ over a field of characteristic $p$, triangularization
over an algebraic closure yields
$\Tr(M^{p^a})=\sum_i\lambda_i^{p^a}=(\sum_i\lambda_i)^{p^a}$.
Reduce the orbit contribution modulo $\pi$ after division to obtain
\eqref{eq:word-leading}.
\end{proof}

For nilpotent residual classes, we shall also need the characteristic
polynomials of all lifts with a prescribed Jordan partition.

\begin{proposition}\label{prop:coefficients}
Let $\lambda=(\lambda_1\ge\cdots\ge\lambda_t)$ be a partition of $m$,
and put $c_i=\min\{j:\lambda_1+\cdots+\lambda_j\ge i\}$.
The characteristic polynomials of integral lifts of $J_\lambda$ are
exactly
\begin{equation}\label{eq:coeff}
 F_x(T)=T^m-\sum_{i=1}^m\pi^{c_i}x_iT^{m-i},\qquad x_i\in\OO.
\end{equation}
Every coefficient choice is realized by an explicit integral matrix.
\end{proposition}
\begin{proof}
The relation between nilpotent reduction and coefficient valuations also
appears in the Newton and Young polygon conditions of Rybakov
\cite{Rybakov}. We give the coefficient argument and construction here,
without a separability hypothesis.
Write a lift as $J_\lambda+\pi X$. In a term of a principal minor of
size $i$, let $j$ factors come from $\pi X$. Deleting those factors
breaks the permutation cycles into at most $j$ paths along the residual
Jordan chains. These paths meet at most $j$ distinct blocks. Thus
$i\le\lambda_1+\cdots+\lambda_j$, and $j\ge c_i$.
Every term of the $i$th characteristic coefficient is divisible by
$\pi^{c_i}$, proving necessity.

For sufficiency, order the blocks by increasing length. In a matrix
$A_x$, put $1$ on the superdiagonal inside each chain and $\pi$ at a
boundary between chains, put $\pi x_i$ in last-row position $m-i+1$,
and put zero in all other positions. Expansion along the last row gives
\[
 \det(TI_m-A_x)=T^m-\sum_{i=1}^m
 (A_x)_{m,m-i+1}\prod_{j=m-i+1}^{m-1}(A_x)_{j,j+1}\,T^{m-i}.
\]
The last $i-1$ superdiagonal entries contain exactly $c_i-1$ boundaries:
they traverse the largest blocks first, from right to left. Thus the
coefficient is $-\pi^{c_i}x_i$. Modulo $\pi$, the boundary entries
and the last row vanish, leaving exactly the specified Jordan chains.
\end{proof}

\section{Nilpotent residual matrices}\label{sec:exact}
Let $\lambda=(\lambda_1\ge\cdots\ge\lambda_t)$ be the residual Jordan
partition. Throughout this section, write $h=\lambda_1$ and
$t_h=\#\{i:\lambda_i=h\}$.
Set $\lambda_2=0$ when there is only one block.
Recall that $n=f(r-1)$, and put $N=p^n$. Then
\[
 L_{h,e}(n)=\min_{\lceil N/h\rceil\le k\le N}
          \{k+e(n-v_p(k))\}.
\]
The exceptional cases in Theorem~\ref{thm:complete} are as follows.
In each case $s,v$ are nonnegative integers and $\nu=n-s-v$;
all the listed conditions are required.
\begin{enumerate}
\renewcommand{\labelenumi}{\textup{(\Roman{enumi})}}
\item\label{case:I} The field and largest block satisfy
$q=2$, $e=2^v$, and $h=3\cdot2^s$. The remaining conditions are
\[
 \nu\ge6,\qquad \nu\equiv0\pmod2,\qquad t_h=1,\qquad s\le\nu-5.
\]
If $s=\nu-5$, require in addition $\lambda_2\le h-2$.
\item\label{case:II} The field and largest block satisfy
$q\in\{2,4\}$, $e=3\cdot2^v$, $\eta=1$, and $h=5\cdot2^s$.
The remaining conditions are
\[
 \nu\ge12,\qquad \nu\equiv0\pmod4,\qquad s\le\nu-9,
 \qquad
 \begin{cases}t_h\le2,&q=2,\\t_h=1,&q=4.\end{cases}
\]
If $s=\nu-9$, require in addition $\lambda_2\le h-2$.
\item\label{case:III} The field and largest block satisfy
$q=p$, $p\equiv1\pmod4$, $e=(p-1)p^v/2$, $\eta^2=-1$, and
$h=(p+1)p^s$. The remaining conditions are
\[
 \nu\ge5,\qquad \nu\equiv1\pmod2,\qquad t_h=1,\qquad s\le\nu-4.
\]
\end{enumerate}
Lemma~\ref{lem:words} gives a lower bound in terms of the largest
Jordan block.
\begin{lemma}\label{lem:nilbound}
Every lift of a nilpotent residual type with largest block $h$ satisfies
\begin{equation}\label{eq:nilbound}
 \vp(\Tr(A^{p^n}))\ge L_{h,e}(n)
 =en+\min_{0\le a\le n}
       \left\{p^a\left\lceil\frac{p^{n-a}}h\right\rceil-ea\right\}.
\end{equation}
For fixed $h$ and $e$, the function $n\mapsto L_{h,e}(n)$ is strictly increasing. In particular
$E_r\ge L_{h,e}(n)\ge\lceil p^n/h\rceil$.
For the standard Jordan lift $B$, the coefficients of
$\Tr((B+\pi X)^{p^n})$, as a polynomial in the entries of $X$,
are all divisible by $\pi^{L_{h,e}(n)}$.
\end{lemma}
\begin{proof}
Here $B^h=0$ and $\Tr(B^j)=0$ for every $j\ge1$.
In a word with $k$ occurrences of $X$, a cyclic gap containing
$h$ successive copies of $B$ gives zero trace. A word with nonzero
trace therefore has $p^n-k\le k(h-1)$, or $k\ge\lceil p^n/h\rceil$.
Lemma~\ref{lem:words} makes every remaining coefficient divisible
by $\pi^kp^{n-v_p(k)}$. This proves the coefficient assertion and
the first lower bound.

For the second expression in \eqref{eq:nilbound}, set
$k_a=p^a\lceil p^{n-a}/h\rceil$. This is an allowed candidate
and $v_p(k_a)\ge a$, so
$L_{h,e}(n)\le k_a+e(n-a)$ for every $a$. Conversely, any allowed
$k$ with $v_p(k)=a$ satisfies $k\ge k_a$, and hence
$k+e(n-a)\ge k_a+e(n-a)$. Taking the two minima proves equality.

For strict increase, consider a candidate $k$ at $n+1$.
If $k=pt$, then $t$ is allowed at $n$, and
\[
 k+e(n+1-v_p(k))
 =t+e(n-v_p(t))+(p-1)t>L_{h,e}(n).
\]
If $p\nmid k$, its value is at least
$\lceil p^n/h\rceil+e(n+1)>L_{h,e}(n)$; here we used
$L_{h,e}(n)\le\lceil p^n/h\rceil+en$.
Thus $L_{h,e}(n+f)>L_{h,e}(n)$, which gives the trace-difference
bound. Finally $v_p(k)\le n$ makes the second summand in the
defining minimum nonnegative.
\end{proof}

For a monic polynomial $F$, write $S_j(F)$ for the sum of the $j$th powers
of its roots, with multiplicity. Newton's identities make $S_j(F)$ an
integral polynomial in the coefficients of $F$. In particular it is
$\Tr(A^j)$ if $F$ is the characteristic polynomial of $A$.

\begin{proof}[Proof of Corollary~\ref{cor:powerblock}]
Suppose $h=p^s$ with $s\ge0$.
When $n\le s$, the least candidate with $v_p(k)=a$ is $k=p^a$.
When $n>s$, the same holds for $a\ge n-s$. An index $a<n-s$ requires
$k>p^{n-s}$ and gives $k+e(n-a)>p^{n-s}+es$.
Thus minimizing the lower bound gives
$L_{h,e}(n)=\min_{0\le j\le\min(s,n)}(ej+p^{n-j})$.
For each allowed $j$, take the companion matrix of
$T^{h-p^j}(T^{p^j}-\pi)$ on one largest block and standard nilpotent
lifts on all other blocks. Its reduction has the prescribed type, and
\[
 \Sps_{p^n}\bigl(T^{h-p^j}(T^{p^j}-\pi)\bigr)
       =p^j\pi^{p^{n-j}}.
\]
Choose $j$ minimizing the valuation. The strict increase of $L_{h,e}$
ensures that $\Tr(A^{q^r})$ has strictly larger valuation.
\end{proof}

For the remaining block lengths, we use the coefficient parameters
$F_x$ and $c_i$ from Proposition~\ref{prop:coefficients}.
Put $\delta_i=hc_i-i\ge0$;
this is the difference between $i$ and the maximum total size $hc_i$
of $c_i$ blocks of length at most $h$.
For a formal power series $H$, $[z^j]H$ denotes its coefficient of $z^j$.
Introduce a formal variable $z$ and the monic polynomial
\begin{equation}\label{eq:aux}
 Q_x(T,z)=T^t-\sum_{j=1}^t B_j(z)T^{t-j},\qquad
 B_j(z)=\sum_{c_i=j}x_i z^{\delta_i}.
\end{equation}
For a formal variable $Z$, the identity $-\log(1-\sum_j B_j Z^j)=\sum_{k\ge1}S_k(Q_x)Z^k/k$
gives the expansion
\begin{equation}\label{eq:costidentity}
 S_N(F_x)=\sum_{\lceil N/h\rceil\le k\le N}
 \pi^k\frac Nk\,[z^{hk-N}]S_k(Q_x).
\end{equation}
To verify \eqref{eq:costidentity}, consider a monomial
$x^\gamma=\prod_i x_i^{\gamma_i}$ with $\sum_i i\gamma_i=N$.
Its weight is $k=\sum_i c_i\gamma_i$, and its exponent of $z$ is
$\sum_i\delta_i\gamma_i=hk-N$. After factoring out $\pi^k$,
its coefficient on both sides is
$N(|\gamma|-1)!/\prod_i\gamma_i!$, where
$|\gamma|=\sum_i\gamma_i$ is its total degree in the variables $x_i$.
Since $S_k(Q_x)$ has integral coefficients, the $k$th summand has valuation
at least $k+e(n-v_p(k))$. This expansion retains all coefficients in
\eqref{eq:coeff}, including those arising from entries between different
residual blocks.

Specializing to a single largest block gives the following necessary
conditions for $E_r>L_{h,e}(n)$.

\begin{proposition}\label{prop:possiblecases}
If $E_r>L_{h,e}(n)$, then the field conditions, the stated form of
$h$, and the lower bound and congruence condition on $\nu$ in one of
\textup{(I)--(III)} hold. The restrictions on $t_h,s,\lambda_2$
will be determined below.
\end{proposition}
\begin{proof}
Corollary~\ref{cor:powerblock} excludes a power of $p$.
Write $h=p^s d$ with $d>1$ and $p\nmid d$, and put $L=L_{h,e}(n)$.
Every lift of a single nilpotent block of length $h$ has characteristic
polynomial $F(T)=T^h-\pi\sum_{i=1}^h x_iT^{h-i}$, and every choice of
$x_i\in\OO$ is realized by its companion matrix. Here $c_i=1$,
so the weight of a monomial equals its total degree in the coefficient
parameters. The power-sum identity is
\begin{equation}\label{eq:newton}
 \Sps_N(F)=\sum_{\sum i\gamma_i=N}
 c_\gamma\pi^k x^\gamma,\qquad
 k=\sum\gamma_i,\qquad c_\gamma=\frac Nk\binom{k}{\gamma_1,\ldots,\gamma_h}.
\end{equation}
Here $x^\gamma=\prod_i x_i^{\gamma_i}$, and $\binom{k}{\gamma_1,\ldots,\gamma_h}$
is the multinomial coefficient. The coefficients $c_\gamma$ are integers.
Moreover
\begin{equation}\label{eq:content}
 v_p(c_\gamma)\ge n-v_p(k),\qquad
 v_p(c_\gamma)\ge n-\min_{\gamma_i>0}v_p(\gamma_i).
\end{equation}
The second inequality follows from
$c_\gamma\gamma_i/N=\binom{k-1}{\gamma_1,\ldots,\gamma_i-1,\ldots,\gamma_h}$.
If $v_p(k)=a$ and the first inequality is an equality, the multinomial
coefficient is a $p$-adic unit. Addition of the $\gamma_i$ in base $p$ then
has no carries, so every $\gamma_i$ is divisible by $p^a$.

By Lemma~\ref{lem:nilbound},
\begin{equation}\label{eq:L}
 L=en+\min_{0\le a\le n}\{k_a-ea\},\qquad
 k_a=p^a\ell_a,\quad \ell_a=\lceil p^{n-a}/h\rceil,
 \quad d_a=h\ell_a-p^{n-a}.
\end{equation}
Every minimizing $a$ has $p\nmid\ell_a$: otherwise increasing $a$ would
give the same $k_a$ and a smaller value. Choose the largest $d_a$ among the
minimizing indices, and call it $d_*$. For
$F=T^h-\pi xT^{d_*}-\pi y$, equations \eqref{eq:newton}--\eqref{eq:content} give
\begin{equation}\label{eq:leading}
 \overline{\pi^{-L}\Sps_N(F)}=
 \sum_{\substack{a\text{ minimizes }\eqref{eq:L}\\d_a=d_*}}
 \eta^{n-a}x^{p^a}y^{k_a-p^a}.
\end{equation}
Indeed, if the multiplicity of $x$ is $p^a t$, the degree constraint gives
$d_* t=d_a$. Maximality of $d_*$ forces $t=1$ and $d_a=d_*$.
The normalized coefficient is $1$ modulo $p$, since
$\binom{p^a\ell_a}{p^a}\equiv\ell_a\pmod p$.

A nonzero value of \eqref{eq:leading} gives $E_r=L$ for the full residual
type: use this companion matrix on one largest block and standard nilpotent
lifts on the remaining blocks. Hence $E_r>L$ forces this function to vanish.
If a term in \eqref{eq:leading} has $\ell_a=1$, setting $y=0$ isolates it.
Otherwise all its indices satisfy $n-a\ge s$. Put
$\tau=\operatorname{ord}_d(p)$, the multiplicative order of $p$ modulo $d$. Indices with the same $d_a$ differ by multiples
of $\tau$, and their values in \eqref{eq:L} are
$en+N/h+(d_*/h)p^a-ea$. Strict convexity along this arithmetic progression
shows that at most two can minimize, and two minimizers are adjacent in
the progression. Thus their indices are $a,a+\tau$. A one-term sum in
\eqref{eq:leading} is nonzero somewhere, for example at $x=y=1$.
For two minimizing indices satisfying $d_a=d_*$, the normalized
multinomial coefficients are units and the weights differ by $e\tau$.
The resulting two-term sum vanishes on all of $\FF_q^2$ only if
\begin{equation}\label{eq:conditions}
 f\mid \tau,\qquad q-1\mid e\tau,\qquad \eta^\tau=-1.
\end{equation}
These conditions follow by varying $x,y\ne0$ independently, using
$k_{a+\tau}-k_a=e\tau$.

We next bound $\tau$. Since
$v_p(k_a)=a$ and $p^{a+\tau}\mid k_{a+\tau}$,
\[
 a=v_p(e\tau),\qquad M=e\tau/p^a,\qquad p\nmid M,
 \qquad \frac{d_*}{h}=\frac{M}{p^\tau-1}.
\]
For $\tau>1$, let $r_0\in\{1,\ldots,p-1\}$ be the residue of $M$
modulo $p$. The next multiple of $p^{a+1}$ gives
\[
 k_{a+1}=k_a+p^ar_0,\qquad
 (k_{a+1}-e(a+1))-(k_a-ea)=p^a(r_0-M/\tau)\ge0.
\]
Equality would make $a+1$ a minimizing index with
$d_{a+1}=(d_*+hr_0)/p>d_*$, since $p^\tau-1>\tau(p-1)$.
This contradicts the choice of $d_*$. Thus $M<\tau r_0\le\tau(p-1)$. For $\tau=1$, $d_*<h$ gives $M<p-1$, contradicting
$p-1\mid M$ from \eqref{eq:conditions}. Write $M=(p-1)k$, $1\le k<\tau$,
and $R=1+p+\cdots+p^{\tau-1}$. Since $\gcd(h,d_*)=p^s$,
\[
 g=\gcd(k,R),\qquad d=R/g,\qquad d_*/p^s=k/g.
\]
For some $0\le j<\tau$, the congruence $p^j\equiv-d_*/p^s\pmod d$ therefore gives
$gp^j+k=\ell R$ with $1\le\ell\le g\le k$. Multiplication by $p-1$
shows $p^j\mid (p-1)k+\ell$, whence $p^j\le pk$ and
\begin{equation}\label{eq:Tbound}
 R\le(\tau-1)\{p(\tau-1)+1\}.
\end{equation}
For $p\ge3$, inequality \eqref{eq:Tbound} fails at $\tau=4$
(the smallest case is $40>30$), and its left side then grows faster
than the quadratic right side. Thus $\tau\le3$.
For $p=2$ it fails at $\tau=6$ ($63>55$) and thereafter, so $\tau\le5$.
For $\tau=3$, $k=1,2$ and $g=1$, but $R-p^j\ge p+1>k$, impossible.
For $p=2,\tau=5$, $R=31$, $k\le4$, and $g=1$, giving the same contradiction.
For $\tau=2$ one obtains $d=p+1$, $M=p-1$, and $d_*=p^s$.
For $p=2,\tau=4$, $k$ is odd; $k=1$ is impossible and $k=3$ gives
$d=5$, $M=3$, and $d_*=p^s$. Conditions \eqref{eq:conditions} now give
the following alternatives. For $\tau=2$, the equality
$e\tau=p^a(p-1)$ gives $e=(p-1)p^a/2$ and $d=p+1$.
If $p=2$, integrality forces $a\ge1$; putting $v=a-1$ gives
case~\textup{(I)}, and $q-1\mid2^a$ forces $q=2$.
If $p$ is odd, $f\mid2$ and $q-1\mid(p-1)p^a$ force $f=1$.
The condition $\eta^2=-1$ then forces $p\equiv1\pmod4$, giving case~\textup{(III)}
with $v=a$.
For $p=2$, $\tau=4$, one has $4e=3\cdot2^a$, so $a\ge2$ and
$e=3\cdot2^{a-2}$. Putting $v=a-2$ gives case~\textup{(II)}.
Here $f\mid4$ and $q-1\mid3\cdot2^a$ leave $q=2,4$;
in either field $\eta^4=1$ forces $\eta=1$.

Finally, $\ell_{a+\tau}>1$ gives the lower bounds and congruence
conditions on $\nu$ in \textup{(I)--(III)}.
For example, in~\textup{(I)}, $a=v+1$, $\tau=2$, and
$\ell_{a+\tau}=(2^{n-s-v-3}+1)/3>1$, so $n-s-v\ge6$ is even.
The corresponding exponents are $n-s-v-6\ge6$ in~\textup{(II)} and
$n-s-v-2\ge3$ in~\textup{(III)}, with the indicated congruences.
\end{proof}

For the passage from power sums to trace differences we use
\begin{equation}\label{eq:separation}
 L_{h,e}(n+1)-L_{h,e}(n)\ge\lfloor(p-1)p^n/h\rfloor.
\end{equation}
For a candidate divisible by $p$, divide it by $p$ and compare with $L_{h,e}(n)$;
the increase is at least $(p-1)\lceil p^n/h\rceil$. A candidate prime to $p$
has value at least $\lceil p^{n+1}/h\rceil+e(n+1)$, whereas
$L_{h,e}(n)\le\lceil p^n/h\rceil+en$. This proves \eqref{eq:separation}.
Under the conditions on $\nu$ in \textup{(I)--(III)}, this lower bound
exceeds $2e$ in~\textup{(I)}, $e/3$ in~\textup{(II)}, and $2e/(p-1)$ in~\textup{(III)}. Hence $\Tr(A^{p^{n+f}})$ has valuation strictly greater than
$L_{h,e}(n)$ plus any correction in \eqref{eq:answer}.

To determine the restrictions on the remaining Jordan blocks, we need
the terms of least valuation in \eqref{eq:costidentity}.

\Needspace{22\baselineskip}
\begin{lemma}\label{lem:weights}
Assume the field conditions, the stated form of $h$, and the lower
bound and congruence condition on $\nu$ from one of
\textup{(I)--(III)}, without imposing the restrictions on
$t_h,s,\lambda_2$. Put $\nu=n-s-v$, $w=p^v$, $H=p^s$,
$d=h/p^s$, and $L=L_{h,e}(n)$.
The weights $k$ in \eqref{eq:costidentity} for which
$k+e(n-v_p(k))$ is below the indicated precision are exactly those
listed below. In each row, $k=p^a\ell$ and $hk-N=p^aH$.
\begin{center}\small
\begin{tabular}{@{}ccccc@{}}\toprule
Case&$a$&$\ell$&valuation minus $L$&precision\\\midrule
\textup{(I)}, $e=w\ge2$&$v+1$&$(2^{\nu-1}+1)/3$&$0$&$L+2e$\\
 &$v+3$&$(2^{\nu-3}+1)/3$&$0$&\\
 &$v-1$&$(2^{\nu+1}+1)/3$&$3e/2$&\\[2pt]
\textup{(II)}, $e=3w$&$v+2$&$(2^{\nu-2}+1)/5$&$0$&$L+w$\\
 &$v+6$&$(2^{\nu-6}+1)/5$&$0$&\\[2pt]
\textup{(III)}, $e=(p-1)w/2$&$v$&$(p^\nu+1)/(p+1)$&$0$&$L+w$\\
 &$v+2$&$(p^{\nu-2}+1)/(p+1)$&$0$&\\\bottomrule
\end{tabular}
\end{center}
In \textup{(I)} with $e=1$, only the first two listed weights contribute
modulo $\pi^{L+1}$. The values of $L$ are
\begin{equation}\label{eq:Lvalues}
 L=\begin{cases}
 es+e(2^\nu/3+\nu-1/3),&\textup{(I)},\\
 es+w(2^\nu/5+3\nu-26/5),&\textup{(II)},\\
 es+w(p^\nu/(p+1)+(p-1)\nu/2+1/(p+1)),&\textup{(III)}.
 \end{cases}
\end{equation}
\end{lemma}
\begin{proof}
Write $L_*$ for the asserted value in \eqref{eq:Lvalues}. For a fixed
$a=v+j\ge0$, a weight with $v_p(k)=a$ has the form $k=p^a\ell$,
where $p\nmid\ell$ and $\ell\ge\lceil p^{\nu-j}/d\rceil$.
For $j\le\nu$, write
$\lceil p^{\nu-j}/d\rceil=(p^{\nu-j}+\epsilon_j)/d$.
The difference between $k+e(n-a)$ and $L_*$, divided by $w$, is
at least
\[
 \begin{cases}
 2^j\epsilon_j/3-j+1/3,&\textup{(I)},\\
 2^j\epsilon_j/5-3j+26/5,&\textup{(II)},\\
 p^j\epsilon_j/(p+1)-(p-1)j/2-1/(p+1),&\textup{(III)}.
 \end{cases}
\]
Equality holds for the least multiple of $p^a$; it is an admissible
weight only if its quotient by $p^a$ is prime to $p$.

In~\textup{(I)}, $\epsilon_j=2$ for even $j$ and $1$ for odd $j$.
For $j\le-2$ the displayed bound exceeds $2$. Its values at
$j=-1,0,1,2,3,4,5$ are respectively
$3/2,1,0,1,0,7,6$. At $j=0,2$ the least quotient is even;
replacing it by the next odd integer raises these bounds to $2$ and $5$.
Increasing $j$ by two changes the displayed expression by
$2^j\epsilon_j-2$, which is positive for $j\ge2$.
Thus only $j=-1,1,3$ can contribute below $L_*+2w$.
At each of these indices the least quotient is odd, and the next
admissible weight raises the bound by $2^{j+1}w$; it is already
above the cutoff. The index $j=-1$ occurs only when $v\ge1$.
If $v=0$, the same calculation leaves only $j=1,3$ below $L_*+1$.

In~\textup{(II)}, the values of $\epsilon_j$ for
$j\equiv0,1,2,3\pmod4$ are $4,2,1,3$.
Negative $j$ give a bound greater than $1$. At $j=0,\ldots,7$
the displayed expression equals $6,3,0,1,6,3,0,61$.
Its change from $j$ to $j+4$ is $3\cdot2^j\epsilon_j-12$,
which is positive for $j\ge4$.
Hence only $j=2,6$ lie below $L_*+w$. The least quotients are odd,
and subsequent admissible weights add at least $2^{j+1}w$.
The same calculation shows that at the endpoint $L_*+w$ there
is exactly one further weight, with $j=3$ and
$\ell=(2^{\nu-3}+3)/5$.

In~\textup{(III)}, $\epsilon_j=1$ for even $j$ and $p$ for odd $j$.
Since $p\ge5$, negative $j$ give a bound greater than $1$.
The values at $j=0,1,2,3,4$ are
$0,(p-1)/2,0,(p-1)(p^2-1/2),(p-1)(p^2-1)$.
The change from $j$ to $j+2$ is
$(p-1)(p^j\epsilon_j-1)>0$ for $j\ge1$.
Only $j=0,2$ therefore lie below $L_*+w$; their least quotients
are prime to $p$, and any later weight adds at least $p^jw\ge w$.

For $j\ge\nu$, the least admissible weight is $p^{v+j}$.
Its bound increases by $(p-1)p^{v+j}-e>0$ at each step,
and the value at $j=\nu$ is already above the relevant cutoff.
The listed weights all lie between $\lceil N/h\rceil$ and $N$.
Since their minimum bound is $L_*$, we have $L=L_*$, and the table follows.
\end{proof}

\Needspace{7\baselineskip}
The next proposition proves that the partition conditions in
\textup{(I)--(III)} are necessary.
\begin{proposition}\label{prop:repeated}
Under the assumptions of Lemma~\ref{lem:weights}, we have
$E_r=L_{h,e}(n)$ if at least one of the following holds:
\begin{enumerate}
\item $t_h\ge2$ in \textup{(I)} or \textup{(III)}; or, in
\textup{(II)}, $t_h\ge2$ when $q=4$ and $t_h\ge3$ when $q=2$.
\item $s\ge\nu-4$ in \textup{(I)}, $s\ge\nu-8$ in
\textup{(II)}, or $s\ge\nu-3$ in \textup{(III)}.
\item $\lambda_2\ge h-1$ and either $s=\nu-5$ in
\textup{(I)} or $s=\nu-9$ in \textup{(II)}.
\end{enumerate}
\end{proposition}
\begin{proof}
For the first assertion, let $j$ be the required number of largest
blocks ($j=2$ or $3$), and take the unramified extension
$L'/K$ of degree $j$. Restriction of scalars of a companion lift of one
block of length $h$ over $\OO_{L'}$ reduces to $j$ copies of $J_h(0)$.
Use the polynomial in \eqref{eq:leading}, with $x,y\in\OO_{L'}$.
Here $d=h/p^s$ and $\tau=\operatorname{ord}_d(p)$, as in the
classification above. Let $a$ and $b=a+\tau$ be the two minimizing
indices, and put $c_a=k_a-p^a$ and $c_b=k_b-p^b$.
For fixed nonzero $y\in\FF_{q^j}$, the normalized residue of the trace
over $K$ is the function
\[
 X\longmapsto\Tr_{\FF_{q^j}/\FF_q}
       (\eta^{n-a}Xy^{c_a}+\eta^{n-b}X^{p^\tau}y^{c_b}),\qquad X=x^{p^a}.
\]
Since $f\mid \tau$, moving Frobenius through the trace makes this a linear
functional in $X$. It is nonzero for some $y$ unless
\[
 q^j-1\mid p^\tau k_a-k_b
       =p^{n-s}\frac{p^\tau-1}{d}.
\]
Here we used $\eta^\tau=-1$ and the nondegenerate finite-field trace pairing.
The prime-to-$p$ factor on the right is $1,3,p-1$ in \textup{(I)}, \textup{(II)}, and \textup{(III)}, respectively.
The stated choices of $j$ violate this divisibility. Choose $y$ and then
$X$ for which the trace is nonzero, and lift them integrally.
Add standard nilpotent lifts on any remaining blocks. This produces
$\vp(\Tr(A^N))=L_{h,e}(n)$, whereas $\Tr(A^{q^r})$ has larger valuation.

For the second assertion, put $H=p^s$ and let $b$ be the larger
of the two minimizing indices in
Lemma~\ref{lem:weights}, and write $k_b=p^b\ell_b$ with $p\nmid\ell_b$.
If $H\ge\ell_b$, then
\[
 (h-1)k_b-N=p^b(H-\ell_b)\ge0,
 \qquad L_{h-1,e}(n)=L_{h,e}(n).
\]
Here $H\ge\ell_b$ means $s\ge \nu-4$ in~\textup{(I)}, $s\ge \nu-8$ in~\textup{(II)},
or $s\ge \nu-3$ in~\textup{(III)}. In particular $h-1$ is prime to $p$ and is
outside the exceptional block lengths for the same field parameters.
Proposition~\ref{prop:possiblecases} therefore gives its exact exponent $L$.
If $G$ is an attaining
characteristic polynomial of degree $h-1$, the ordinary companion matrix
of $TG$ has residual type $J_h(0)$ and the same positive power sums.
Adjoin standard nilpotent blocks for the remaining parts of $\lambda$.
This gives $E_r=L$ whenever $H\ge\ell_b$. The complementary inequality
$H<\ell_b$ is exactly
\begin{equation}\label{eq:srestrictions}
 s\le \nu-5\quad\textup{(I)},\qquad
 s\le \nu-9\quad\textup{(II)},\qquad
 s\le \nu-4\quad\textup{(III)}.
\end{equation}
The three comparisons follow from
\[
 \begin{aligned}
 2^{\nu-5}&<\frac{2^{\nu-3}+1}{3}<2^{\nu-4},
 &&\text{in \textup{(I)}},\\
 2^{\nu-9}&<\frac{2^{\nu-6}+1}{5}<2^{\nu-8},
 &&\text{in \textup{(II)}},\\
 p^{\nu-4}&<\frac{p^{\nu-2}+1}{p+1}<p^{\nu-3},
 &&\text{in \textup{(III)}}.
 \end{aligned}
\]

Finally, at $s=\nu-5$ in~\textup{(I)}, or $s=\nu-9$ in~\textup{(II)}, suppose
$\lambda_2\ge h-1$. Both coefficient indices $h-1$ and $2h-1$ then occur,
with $c_{h-1}=1$, $c_{2h-1}=2$, and
$\delta_{h-1}=\delta_{2h-1}=1$. Choose only these two coefficients
nonzero, with $x_{h-1}=x_{2h-1}=1$ in \eqref{eq:coeff}.
A term at weight $k_b=p^b\ell_b$ must have
\begin{equation}\label{eq:boundaryterm}
 \gamma_{h-1}=p^b(2H-\ell_b),\qquad
 \gamma_{2h-1}=p^b(\ell_b-H).
\end{equation}
Here $H<\ell_b<2H$, and $j=\ell_b-H$ is odd:
$j=(H+1)/3$ in~\textup{(I)} and $j=(3H+1)/5$ in~\textup{(II)}.
The coefficient is
\[
 \frac{2^n}{2^{b+s}}\binom{2^{b+s}}{2^b j},\qquad
 v_2\binom{2^{b+s}}{2^b j}=s,
\]
so its valuation is exactly $n-b$. At the other minimum weight one has
$\ell_a>2H$, making \eqref{eq:boundaryterm} impossible. Thus there is
exactly one lowest term, of valuation $L$. This proves $E_r=L$ in the
boundary cases as well, for a lift of the entire prescribed partition.
The strict increase of $L$ gives $\vp(\Tr(A^{q^r}))>L$.
\end{proof}

\Needspace{7\baselineskip}
The lower bounds in the remaining cases require comparing power sums
for different choices of the coefficient parameters.

\begin{lemma}\label{lem:precision}
If $a=v_p(k)$ and $x_i'\equiv x_i\pmod\pi$, then coefficientwise
\[
 S_k(Q_{x'})\equiv S_k(Q_x)\pmod{\pi^{C_a}},\qquad
 C_a=\min_{0\le j\le a}\{p^j+e(a-j)\}.
\]
For binary coefficients $x_i\in\{0,1\}$ one also has, for $a\ge1$ and
odd $\ell$,
\begin{equation}\label{eq:modfour}
 S_{2^a\ell}(Q_x(T,z))\equiv
 S_{2\ell}(Q_x(T,z^{2^{a-1}}))\pmod4.
\end{equation}
Over a field of characteristic $p$, $S_{p^a\ell}(Q_x)=S_\ell(Q_x)^{p^a}$.
\end{lemma}
\begin{proof}
Write $B(Z,z)=\sum B_j(z)Z^j$ and replace $B$ by $B+\pi W$.
The coefficient of the $j$th power of a perturbation parameter in $S_k$
is
\[
 \frac kj[Z^k]\frac{W^j}{(1-B)^j}.
\]
It is an integral polynomial, because $S_k$ is one; its coefficients also
have $p$-adic valuation at least $a-v_p(j)$. After multiplication by
$\pi^j$, the valuation is at least
$j+e\max\{0,a-v_p(j)\}\ge C_a$. Summing proves the first assertion.
For \eqref{eq:modfour}, use the companion matrix of $Q_x$. Its entries
are sums of monomials in $z$ with coefficients $0$ or $1$. Group closed
walks by their primitive cyclic periods. In
$S_{2^a\ell}(z)-S_{2^{a-1}\ell}(z^2)$, the periods dividing
$2^{a-1}\ell$ cancel; every remaining period is divisible by $2^a$.
The difference is therefore divisible by $2^a$. Iterate for $a\ge2$;
when $a=1$ the assertion is an identity. The characteristic-$p$ assertion
follows by taking powers of the roots, or by the same closed-walk expansion.
\end{proof}

\begin{proof}[Completion of the proof of Theorem~\ref{thm:complete}]
Propositions~\ref{prop:possiblecases} and \ref{prop:repeated} give
$E_r=L_{h,e}(n)$ outside \textup{(I)--(III)}. Assume henceforth
all the conditions in one of these three cases.
Recall from \eqref{eq:coeff}--\eqref{eq:aux} that
$c_i=\min\{j:\lambda_1+\cdots+\lambda_j\ge i\}$,
$\delta_i=hc_i-i$, and $Q_x(T,z)=T^t-\sum B_j(z)T^{t-j}$, where
$B_j(z)=\sum_{c_i=j}x_i z^{\delta_i}$. Write $\nu=n-s-v$, $H=p^s$, $w=p^v$, and $L=L_{h,e}(n)$.
Let $a<b$ be the two minimizing indices in Lemma~\ref{lem:weights},
and write $k_a=p^a\ell_a$ and $k_b=p^b\ell_b$.
We compare the coefficients at these two weights.
In \textup{(I)}, the comparison is made modulo $4$ and also involves
the third weight in Lemma~\ref{lem:weights}.

The equality $\delta_i=0$ holds precisely when $i=jh$ with
$1\le j\le t_h$.
Reduce $x_i$ to the residue field and let $t_0\le t_h$ be the largest $j$
with $x_{jh}\ne0$; take $t_0=0$ when all these coefficients vanish.
At $z=0$, \eqref{eq:aux} is $T^{t-t_0}$ times a monic polynomial of degree
$t_0$ with nonzero constant coefficient. Hensel factorization in
$\FF_q[[z]]$ separates these two coprime factors. Denote by $Q_1$ the factor whose roots have $z$-adic valuation zero,
and by $Q_0$ the factor whose roots have positive $z$-adic valuation.
These valuations are taken in an algebraic closure of $\FF_q((z))$.
For $t_0=0$, we have $Q_1=1$.

The coefficient constraints imply
\begin{equation}\label{eq:delta}
 \delta_i\ge c_i-t_h\quad(c_i>t_h),
 \qquad \delta_i\ge1\text{ whenever }\delta_i\ne0.
\end{equation}
Indeed $i\le ht_h+(c_i-t_h)(h-1)$ for $c_i>t_h$.
Use the Newton polygon formed from the points
$(j,\operatorname{ord}_z B_j)$ and $(0,0)$, where
$\operatorname{ord}_z$ denotes the order of vanishing at $z=0$.
After its endpoint $(t_0,0)$ on the horizontal axis, every slope is at
least $1/(t_h+1-t_0)$: for $t_0<c_i\le t_h$ use
$1/(c_i-t_0)$, and for $c_i>t_h$ use
$(c_i-t_h)/(c_i-t_0)\ge1/(t_h+1-t_0)$.
Consequently
\begin{equation}\label{eq:nilaux}
 \operatorname{ord}_z S_\ell(Q_0)\ge
 \left\lceil\frac{\ell}{t_h+1-t_0}\right\rceil.
\end{equation}
If $t_h=1$ and $\lambda_2\le h-2$, the stronger inequality
$\delta_i\ge2(c_i-1)$ for $c_i\ge2$ makes every positive slope at least $1$,
even when $t_0=0$.

In~\textup{(I)}, $t_h=1$. Away from its boundary, $\ell_b>2H$; at the boundary the
stronger slope bound applies. In~\textup{(II)}, away from its boundary $t_h\le2$ and
$\ell_b>3H$; at its boundary $\lambda_2\le h-2$, so $t_h=1$ and the stronger
bound applies. In~\textup{(III)}, $t_h=1$ and $\ell_b>2H$.
Thus, in all remaining cases,
\begin{equation}\label{eq:dropnil}
 [z^H]S_\ell(Q_x)=[z^H]S_\ell(Q_1)
 \quad\text{for }\ell=\ell_a,\ell_b.
\end{equation}
For~\textup{(I)} with binary integer coefficients, the same factorization exists in
$\ZZ[[z]]$: when $t_0=1$, the factor $Q_1$ is $T-V(z)$ with $V(0)=1$;
when $t_0=0$, we have $Q_1=1$. Apply the slope bounds over
$\mathbb Q((z))$. They also show that $Q_0$, whose roots have positive $z$-adic valuation,
contributes nothing to the coefficient of $z^{2H}$ in $S_{2\ell}$,
or to the coefficient of $z^H$ in the third power sum. This gives the integral analogue of \eqref{eq:dropnil} used in
case~\textup{(I)}.

Consider cases \textup{(II)} and \textup{(III)}.
By Lemma~\ref{lem:precision}, the two minimum-weight terms depend only on
Teichm\"uller representatives modulo $\pi^w$. Since $e\ge w$, these
representatives modulo $\pi^w$ form a copy of $\FF_q$.
The normalized expression is
\begin{equation}\label{eq:II-III-residue}
 \pi^{-L}S_N(F_x)\equiv
 u^{n-a}\ell_a^{-1}([z^H]S_{\ell_a}(Q_x))^{p^a}
 +u^{n-b}\ell_b^{-1}([z^H]S_{\ell_b}(Q_x))^{p^b}
 \pmod{\pi^w}.
\end{equation}
Its coefficient expressions on the right are interpreted in that copy of
$\FF_q$. The two denominators agree modulo $p$, and $f\mid b-a$.

In~\textup{(II)}, $\ell_a-\ell_b=3\cdot2^{\nu-6}$. If $q=2$, the companion matrix $M(z)$
of $Q_1$ has size at most $2$, with $M(0)\in\mathrm{GL}_{t_0}(\FF_2)$.
Since $M(0)$ has order dividing $6$, write
$M(z)^6=I+zR(z)$ with $R(z)$ a matrix over $\FF_2[[z]]$.
Raising this identity to the $2^{\nu-7}$th power gives
\[
 M(z)^{3\cdot2^{\nu-6}}\equiv I\pmod{z^{2^{\nu-7}}}.
\]
The exponent $2^{\nu-7}$ exceeds $H$, since $s\le \nu-9$.
If $q=4$, the factor $Q_1$ has degree at most $1$ and its nonzero root at $z=0$
has order dividing $3$, giving the stronger modulus $z^{2^{\nu-6}}$.
Thus \eqref{eq:dropnil} implies
$[z^H]S_{\ell_a}(Q_x)=[z^H]S_{\ell_b}(Q_x)$.
In~\textup{(III)}, the factor $Q_1$ has degree at most $1$, and
$\ell_a-\ell_b=(p-1)p^{\nu-2}$; its nonzero residue lies in $\FF_p^\times$.
The same argument gives equality modulo $z^{p^{\nu-2}}$, again beyond degree $H$.
If $t_0=0$, both coefficients are zero and no matrix is required.
Therefore \eqref{eq:II-III-residue} has respectively the factors $u^4+1$ and
$u^2+1$, proving the lower bounds in \eqref{eq:answer}.

In case~\textup{(I)}, assume first that $e=2^v\ge2$. Denote the three factors in Lemma~\ref{lem:weights} by
$\ell_a,\ell_b,\ell_c$, at indices $a=v+1$, $b=v+3$, $c=v-1$.
Lemma~\ref{lem:precision} permits replacing every $x_i$ by $0$ or $1$
modulo $\pi^{2e}$ in the first two terms, and modulo $\pi^{e/2}$ in the third.
Indeed $C_a\ge2e$, $C_b\ge2e$, and $C_c\ge e/2$.
Equation \eqref{eq:modfour} and its characteristic-two counterpart give
\begin{equation}\label{eq:Iterms}
 \pi^{-L}S_N(F_x)\equiv
 u^{n-a}\alpha_a+u^{n-b}\alpha_b
 +\pi^{3e/2}u^{n-c}\beta\pmod{\pi^{2e}},
\end{equation}
where
\[
 \alpha_j=\ell_j^{-1}[z^{2H}]S_{2\ell_j}(Q_x)\pmod4\quad(j=a,b),\qquad
 \beta=[z^H]S_{\ell_c}(Q_x)\pmod2.
\]
These are integer residue classes. When $Q_1=1$, they are
all zero. Otherwise the integral version of \eqref{eq:dropnil} replaces
these power sums by $V(z)^{2\ell_a}$, $V(z)^{2\ell_b}$, and $V(z)^{\ell_c}$,
where $V\in\ZZ[[z]]$, $V(0)=1$.
The identities
\[
 \ell_a-\ell_b=2^{\nu-3},\qquad \ell_c-\ell_b=5\cdot2^{\nu-3},
 \qquad \nu-3\ge s+2
\]
imply
\begin{equation}\label{eq:Iequal}
 \alpha_a=\alpha_b\pmod4,\qquad \beta=\alpha_b\pmod2.
\end{equation}
For the first congruence, the denominators agree modulo $4$, and every
coefficient through degree $2H=2^{s+1}$ in
$V^{2^{\nu-2}}-1$ is divisible by $4$. Indeed, writing $V=1+zW$ gives
$v_2\binom{2^{\nu-2}}j=\nu-2-v_2(j)\ge\nu-s-3\ge2$ for
$1\le j\le2H$.
For the second, use Frobenius modulo $2$ and the equality of coefficients
through degree $H$ in $V^{\ell_c}$ and $V^{\ell_b}$.
Taking an integer representative $\alpha$ of the common residue in
\eqref{eq:Iequal}, equation \eqref{eq:Iterms} becomes
\[
 \pi^{-L}S_N(F_x)\equiv
 u^{n-b}\alpha\bigl(u^2+1+\pi^{3e/2}u^4\bigr)\pmod{\pi^{2e}}.
\]
Replacing $\beta$ by $\alpha$ changes the expression by a multiple of
$2\pi^{3e/2}$, which vanishes modulo $\pi^{2e}$. This proves the lower bound in case~\textup{(I)}, including even $\alpha$.
When $e=1$, only the first two terms modulo $\pi$ are needed.
Equation \eqref{eq:dropnil} and $\ell_a-\ell_b=2^{\nu-3}>H$ make their
coefficient functions equal over $\FF_2$; their sum is zero. Hence
$\vp(S_N(F_x))\ge L+1$ in that case.

For attainment, consider first $s=0$, so that $h=3$, $5$, or $p+1$
in the respective cases. Write $L=L_{h,e}(n)$.
The polynomial $F_0=T^h-\pi T-\pi$ has auxiliary polynomial
$Q(T,z)=T-(1+z)$. Hence $[z]S_\ell(Q)=\ell$.
In \textup{(II)} and \textup{(III)}, substitution in \eqref{eq:II-III-residue} gives
\[
 \pi^{-L}S_N(F_0)\equiv
 \begin{cases}
 u^{n-v-6}(u^4+1)\pmod{\pi^{p^v}},&\textup{(II)},\\
 u^{n-v-2}(u^2+1)\pmod{\pi^{p^v}},&\textup{(III)}.
 \end{cases}
\]
In~\textup{(I)} with $e\ge2$, one has
$\ell^{-1}[z^2](1+z)^{2\ell}=2\ell-1\equiv1\pmod4$ for odd $\ell$.
Equation \eqref{eq:Iterms} consequently gives
\begin{equation}\label{eq:cubicvalue}
 \pi^{-L}S_N(F_0)\equiv
 u^{n-v-3}\bigl(u^2+1+\pi^{3e/2}u^4\bigr)\pmod{\pi^{2e}}.
\end{equation}
These polynomials attain the required value whenever the correction is
strictly below $2e$ in~\textup{(I)} or $p^v$ in \textup{(II)} and \textup{(III)}. The endpoint values require different coefficient choices.

For case~\textup{(I)}, suppose first that $e=w=2^v\ge2$ and
$\vp(u^2+1+\pi^{3w/2}u^4)\ge2w$. Put
$F_1=T^3-\pi(1+\pi)T-\pi$. For
$F_z=T^3-\pi zT-\pi$, write
\[
 k_a=(N+2w)/3,\qquad k_b=(N+8w)/3,\qquad
 k_c=(N+w/2)/3.
\]
\Needspace{9\baselineskip}
In \eqref{eq:newton}, the multiplicity of $\pi z$ is $j=3k-N$.
The only terms of $S_N(F_z)$ whose coefficients have valuation below
$L+2w$ are as follows:
\begin{center}\small
\begin{tabular}{@{}cccc@{}}\toprule
$k$ & $j$ & valuation minus $L$ & normalized odd coefficient\\\midrule
$k_a$ & $2w$ & $0$ & $1\pmod4$\\
$k_b$ & $8w$ & $0$ & $1\pmod4$\\
$k_c$ & $w/2$ & $3w/2$ & $1\pmod2$\\\bottomrule
\end{tabular}
\end{center}
Indeed, Lemma~\ref{lem:weights} leaves only these three values of
$k$, and $j=3k-N$ fixes the multiplicity. Their normalized coefficients
are $\ell^{-1}\binom{2^a\ell}{2^a}$, where $k=2^a\ell$ and $\ell$ is odd.
For the first two rows $a\ge2$, and this coefficient is $1$ modulo $4$:
in the product
$\prod_{i=1}^{2^a-1}(2^a\ell-i)/(2^a-i)$, the factors with
$v_2(i)\le a-2$ are $1$ modulo $4$, and the remaining factor is
$2\ell-1\equiv1\pmod4$.

The binomial expansion gives
\begin{equation}\label{eq:unitpowers}
 \vp((1+\pi)^{2w}-1)>2w,
 \qquad \vp((1+\pi)^{w/2}-1)=w/2.
\end{equation}
For the first assertion, the only possible terms of valuation $2w$ have
indices $w,2w$; their normalized residues add to $1+1=0$ in $\FF_2$.
For the second assertion, the last binomial term is the unique term of
least valuation. The first assertion also applies with exponent $8w$.
In $S_N(F_1)-S_N(F_0)$, the first two rows therefore have valuation
greater than $L+2w$, and the third has valuation exactly $L+2w$.
Every omitted term has coefficient valuation at least $L+2w$ and gains
at least one further power of $\pi$ in the subtraction. Thus
\[
 \vp\bigl(S_N(F_1)-S_N(F_0)\bigr)=L+2w.
\]
Since both values have valuation at least $L+2w$, at least one attains it.
For $e=1$, the field is $\QQ_2$ and the exceptional cases have even
$n\ge6$. For $F_0=T^3-2T-2$, only the weights $(N+2)/3$ and $(N+8)/3$
contribute modulo $2^{L+2}$. Their normalized coefficients are both
$1$ modulo $4$, so $2^{-L}S_N(F_0)\equiv2\pmod4$ and the correction is $1$.

In case~\textup{(II)}, put $w=2^v$, so $e=3w$.
To attain the endpoint $w$, use $F=T^5-\pi T^3-\pi$.
At weight $k$, the multiplicity of the coefficient of $T^3$ is
$j=(5k-N)/3$. At the two minimum weights, $5k-N$ equals
$4w$ or $64w$, neither divisible by $3$, so both terms are absent.
The proof of Lemma~\ref{lem:weights} leaves exactly one additional
weight whose bound is at most $L+w$, namely
\[
 k=\frac{N+24w}{5}=2^{v+3}\frac{2^{\nu-3}+3}{5},\qquad j=2^{v+3}.
\]
The quotient $\ell=(2^{\nu-3}+3)/5$ is odd. Thus $v_2(k)=v+3$
and $\binom{k}{j}=\binom{2^{v+3}\ell}{2^{v+3}}\equiv\ell\equiv1\pmod2$.
The coefficient $N\binom{k}{j}/k$ therefore has valuation
$e(n-v-3)$, and the corresponding term has valuation $L+w$.
Every other term has larger valuation, so $\vp(S_N(F))=L+w$.

In case~\textup{(III)}, put $w=p^v$ and $e=(p-1)w/2$. Use $F=T^{p+1}-\pi^2T-\pi$.
At weight $k$, the multiplicity of the coefficient of $T$ is
$j=(p+1)k-N$, and the extra factor $\pi$ in this coefficient adds $j$
to the valuation. The weight $k_a=(N+w)/(p+1)$ has $j=w$, a unit
binomial coefficient, and valuation $L+w$.
If $k>k_a$, then $j>w$, so the valuation is greater than $L+w$.
If $k=k_a-t<k_a$ contributes, then $0<t<w/(p+1)$ and
$v_p(k)=v_p(t)<v$. Its valuation minus $L+w$ is at least
\[
 -(p+2)t+e\bigl(v-v_p(t)\bigr)>0,
\]
because $e\ge2w>(p+2)w/(p+1)$. Thus $\vp(S_N(F))=L+w$.

Return to $h=p^s d$, where $d=3,5$, or $p+1$ in the corresponding case.
Apply the preceding constructions at exponent $p^{n-s}$ and let $F_0$
be a polynomial of degree $d$ attaining the required correction.
These power-sum calculations only use the stated conditions on
$\nu=n-s-v$; they do not require $n-s$ to be a multiple of $f$.
The companion matrix of $F_0(T^{p^s})$ reduces to $J_h(0)$, and
\[
 S_{p^n}\bigl(F_0(T^{p^s})\bigr)=p^sS_{p^{n-s}}(F_0),\qquad
 L_{h,e}(n)=es+L_{d,e}(n-s).
\]
Adjoin standard nilpotent lifts of the remaining blocks of $\lambda$.
Their positive power sums vanish. The resulting matrix has the entire
prescribed residual type and attains $L$ plus the correction in
\eqref{eq:answer}. By \eqref{eq:separation}, the trace at $p^{n+f}=qN$
has strictly larger valuation, so the same value is attained by $D_r$.
This proves \eqref{eq:answer} and attainment in all three cases.
\end{proof}

\section{Nonnilpotent residual matrices}\label{sec:nonnil}
The cyclic-permutation estimate gives a criterion for $E_r>C_n$
that applies to any residual matrix. Recall that $C_n=\min_{0\le j\le n}\{p^j+e(n-j)\}$,
where $n=f(r-1)$.
Let $I_n=\{j:0\le j\le n,\ p^j+e(n-j)=C_n\}$.
Recall that $u=p/\pi^e$ and $\eta=\overline u$.
\begin{lemma}\label{thm:first}
Choose $B$ as in Lemma~\ref{lem:periodic}, and write $A=B+\pi X$.
Then $\vp(D_r(A))\ge C_n$, and
\begin{equation}\label{eq:first}
 \pi^{-C_n}D_r(A)\bmod\pi
 =-\sum_{a\in I_n}\eta^{\,n-a}
   \Tr\bigl(\barA^{\,p^{n-a}-1}\overline X\bigr)^{p^a}.
\end{equation}
Consequently, $E_r>C_n$ if and only if the right-hand side is
the zero function on $\Mat_m(\FF_q)$.
\end{lemma}
\begin{proof}
After integral conjugation, we may assume that $A\bmod\pi=B\bmod\pi$.
At exponent $p^n$, Lemma~\ref{lem:words} gives
\[
 \Tr(A^{p^n})-\Tr(B^{p^n})
 =\sum_{t=1}^{p^n}\pi^tH_{n,t}(B,X),\qquad
 \vp(\pi^tH_{n,t}(B,X))\ge t+e(n-v_p(t)).
\]
Write $t=p^ab$, where $p\nmid b$ and $0\le a\le n$. Then
\[
 t+e(n-v_p(t))=p^ab+e(n-a)
 \ge p^a+e(n-a)\ge C_n.
\]
The first inequality is an equality only when $b=1$. Hence a term
can survive after division by $\pi^{C_n}$ and reduction only when
$t=p^a$ and $a\in I_n$. For such an $a$, put
$H_{n,p^a}=p^{n-a}F_a$. The equality
$p^a+e(n-a)=C_n$ and $p=u\pi^e$ imply
\[
 \pi^{-C_n}\pi^{p^a}H_{n,p^a}=u^{n-a}F_a,\qquad
 \overline F_a=\Tr(\barA^{p^{n-a}-1}\overline X)^{p^a}.
\]
The candidates defining $C_{n+1}$ consist of the old candidates
increased by $e$, together with $p^{n+1}$. Since $C_n\le p^n$,
\[
 C_{n+1}=\min(C_n+e,p^{n+1})>C_n.
\]
It follows that $\Tr(A^{p^{n+f}})-\Tr(B^{p^{n+f}})$ has valuation
at least $C_{n+f}>C_n$. Subtract the corresponding difference at exponent
$p^n$ and use
$D_r(B)=0$ to obtain \eqref{eq:first}.
Every matrix $\overline X\in\Mat_m(\FF_q)$ has an integral lift. Thus
a nonzero value of the right side gives a lift with valuation exactly
$C_n$; if it vanishes everywhere, every lift has valuation at least
$C_n+1$. This proves the equivalence.
\end{proof}

\Needspace{13\baselineskip}
The vanishing condition in \eqref{eq:first} can be written as a matrix
identity by examining the indices at which the minimum is attained.

\begin{corollary}\label{cor:improvement}
The set $I_n$ consists of one element or two adjacent elements.
\begin{enumerate}
\item If $I_n=\{a\}$, then $E_r>C_n$ if and only if
$\barA^{p^{n-a}-1}=0$.
\item If $I_n=\{a,a+1\}$, then $e=(p-1)p^a$ and $n\ge a+1$.
Put $b=n-a-1$. If $f>1$, improvement is equivalent to
$\barA^{p^b-1}=0$. If $f=1$, it is equivalent to
\begin{equation}\label{eq:exception}
 \barA^{p^b-1}(\barA^{p-1}+\eta^{-1}I_m)^{p^b}=0.
\end{equation}
\end{enumerate}
Here, as usual, $\barA^0=I_m$.
\end{corollary}
\begin{proof}
Write $F_n(a)=p^a+e(n-a)$ for $0\le a\le n$. Then
\[
 F_n(a+1)-F_n(a)=(p-1)p^a-e.
\]
These differences strictly increase. The minimum is therefore unique
unless one difference is zero, in which case the two adjacent indices
are the only minimizers and $e=(p-1)p^a$.
For a single minimizer, the scalar coefficient and the Frobenius
power in \eqref{eq:first} are invertible. The function vanishes exactly
when $\Tr(\barA^{p^{n-a}-1}X)=0$ for all $X$. The pairing is
nondegenerate: if $M_{ij}\ne0$, then $\Tr(ME_{ji})=M_{ij}$.
This proves the first assertion.

For two minimizers, put $b=n-a-1$ and define the $\FF_q$-linear forms
$L_1(X)=\Tr(\barA^{p^{b+1}-1}X)$ and
$L_2(X)=\Tr(\barA^{p^b-1}X)$. Dividing by
$\eta^{n-a-1}$ and applying the inverse of the $p^a$-power
automorphism reduces the condition to
\[
 cL_1(X)+L_2(X)^p=0\quad\text{for every }X,
 \qquad c^{p^a}=\eta.
\]
If $f>1$, choose $t\in\FF_q$ with $t^p\ne t$. Subtract $t$ times
this identity at $X$ from the identity at $tX$. The result is
$(t^p-t)L_2(X)^p=0$, so $L_2=0$ and then $L_1=0$.
By nondegeneracy, $L_2=0$ means $\barA^{p^b-1}=0$, which also
implies $\barA^{p^{b+1}-1}=0$.
If $f=1$, then $c=\eta$ and $L_2(X)^p=L_2(X)$. Thus the condition
is the matrix identity
\[
 \eta\barA^{p^{b+1}-1}+\barA^{p^b-1}=0
 \iff \barA^{p^b-1}(\barA^{p-1}+\eta^{-1}I_m)^{p^b}=0.
\]
The last equality uses the characteristic-$p$ binomial identity for
commuting matrices and $\eta^{p^b}=\eta$. This is
\eqref{eq:exception}.
\end{proof}

\begin{proof}[Proof of Theorem~\ref{thm:nonmain}]
No nonnegative power of a nonnilpotent matrix is zero.
Corollary~\ref{cor:improvement} therefore shows that $E_r>C_n$
holds precisely under the exceptional condition in
Theorem~\ref{thm:nonmain}. Otherwise, Lemma~\ref{thm:first}
gives $E_r=C_n$ and an integral lift attaining this value.
\end{proof}

\begin{remark}
The residue in \eqref{eq:first} is a function on a finite field. A nonzero
formal polynomial may induce the zero function; over $\FF_p$, for
example, $X^p$ and $X$ induce the same function. This distinction
accounts for the exceptional condition \eqref{eq:exception}.
\end{remark}

For Theorem~\ref{thm:nonexact}, assume $f=1$ and
$e=(p-1)p^a$, and write $n=a+b+1$, where $a,b\ge0$.
The matrix identity \eqref{eq:exception} says precisely that
all nonzero residual eigenvalues satisfy $\alpha^{p-1}=-\eta^{-1}$,
all their Jordan blocks have length at most $p^b$, and
$h_0\le p^b-1$. This follows because $T^{p-1}+\eta^{-1}$ is separable
and relatively prime to $T$.
Set
\begin{equation}\label{eq:delta0}
 \delta_0=p^a-ea,\qquad C_n=en+\delta_0.
\end{equation}
\Needspace{12\baselineskip}
The logarithm on principal units is given by its convergent series;
in a ramified field its values need not be integral.
The next two lemmas estimate the difference between
$\Tr(M^N)-d$ and $N\log\det M$ for matrices with unipotent reduction.

\begin{lemma}\label{lem:error}
Suppose $M\in M_d(\OO)$ has unipotent reduction whose largest Jordan
block has length at most $p^b$. For $N=p^n$, put
\[
 R_N(M)=\Tr(M^N)-d-N\log\det M.
\]
Then $\vp(R_N(M))\ge C_n+e$.
The same bound holds over every unramified extension of $K$.
\end{lemma}
\begin{proof}
Write $X=M-I_d$ and
$\det(TI_d-X)=T^d-\sum_i z_iT^{d-i}$.
If the residual partition of $X$ is $\lambda$, its coefficient constraints
are $z_i=\pi^{c_i}x_i$, where
$c_i=\min\{j:\lambda_1+\cdots+\lambda_j\ge i\}$. In particular
$i\le p^bc_i$, which is the only property of the partition used below.

For a nonzero multi-index $\alpha$, write
\[
 k=\sum_i\alpha_i,\qquad I=\sum_i i\alpha_i,\qquad
 w=\sum_i c_i\alpha_i,\qquad t=\min_{\alpha_i>0}v_p(\alpha_i).
\]
Thus $p^t\mid w$, $w\ge p^t$, and $I\le p^bw$.
Expanding the logarithm in
\[
 \det(I_d-zM)=(1-z)^d\left(1-\sum_i z_i\left(\frac z{1-z}\right)^i\right)
\]
shows that the coefficient of $z_1^{\alpha_1}\cdots z_d^{\alpha_d}$ in
$\Tr(M^N)-d$ is
\[
 \frac Nk\binom{k}{\alpha_1,\ldots,\alpha_d}\binom{N-1}{I-1}.
\]
We take the last binomial coefficient to be zero for $I>N$.
The corresponding coefficient in $N\log\det M$ is
$\frac Nk\binom{k}{\alpha}(-1)^{I-1}$, because
$\det M=1-\sum_i(-1)^i z_i$.
Hence the coefficient in the remainder is
\begin{equation}\label{eq:remcoefficient}
 \frac Nk\binom{k}{\alpha}
       \left\{\binom{N-1}{I-1}-(-1)^{I-1}\right\}.
\end{equation}
Choose $j$ with $v_p(\alpha_j)=t$, and let $\mathbf e_j$ denote the
$j$th coordinate vector. The identity
$k^{-1}\binom{k}{\alpha}=\alpha_j^{-1}\binom{k-1}{\alpha-\mathbf e_j}$
gives valuation at least $n-t$ for its first factor.
For $2\le I\le N$, the identity
\[
 \binom{N-1}{I-1}=(-1)^{I-1}
                 \prod_{j=1}^{I-1}(1-N/j)
\]
gives valuation at least $n-\lfloor\log_p(I-1)\rfloor$ for the expression
in braces. For $I>N$ its valuation is zero; for $I=1$ it vanishes.
Consequently the monomial, after substituting $z_i=\pi^{c_i}x_i$, has
valuation at least
\begin{equation}\label{eq:weightbound}
 en+w-et+e\max\{0,a+2-\lceil\log_p w\rceil\}.
\end{equation}
Indeed $I-1<p^bw\le p^{b+\lceil\log_p w\rceil}$.
The infinite sum coming from the logarithm converges, since $w$ grows
linearly and the denominator losses grow at most logarithmically.

If $w=p^j$, the expression in \eqref{eq:weightbound} is least
when $t=j$, and its part after $en$ is
\[
 p^j-ej+e\max\{0,a+2-j\}\ge\delta_0+e.
\]
For $p>2$ the minimum occurs at $j=a+1$; for $p=2$ it occurs at
$j=a+1,a+2$. This follows by subtracting consecutive values: the differences
are $(p-1)p^j-2e$ for $j<a+2$ and $(p-1)p^j-e$ thereafter.
If $w$ is not a power of $p$, put $j=\lceil\log_p w\rceil$.
Then $t\le j-1$ and $w\ge p^{j-1}+p^t$, so the part after $en$ is at least
\[
 \delta_0+p^{j-1}+e\max\{0,a+2-j\}>\delta_0+e.
\]
Here $p^t-et\ge\delta_0$; if $j\le a+1$ the last term is at least $e$,
and if $j\ge a+2$ one has $p^{j-1}>e$.
All these estimates remain valid over an unramified extension of $K$.
\end{proof}

When $p=2$, cancellation between two terms improves this estimate.

\begin{lemma}\label{lem:binaryerror}
Suppose $p=2$, $f=1$, $a\ge1$, and $M$ satisfies the assumptions of
Lemma~\ref{lem:error}. Then
\begin{equation}\label{eq:binaryerror}
 \vp(R_N(M))\ge C_n+
       \min\{3e/2, e+\vp(u^2+1)\}.
\end{equation}
\end{lemma}
\begin{proof}
Use the coefficient expansion and the notation $k,I,w,t$ from the
proof of Lemma~\ref{lem:error}. For $w$ not a power of $2$, put
$j=\lceil\log_2w\rceil$. The part of the bound
\eqref{eq:weightbound} after $en$ is at least
$\delta_0+2e$: if $j\le a$, the last term is at least $2e$;
if $j=a+1,a+2$, the terms $2^{j-1}+e\max\{0,a+2-j\}$ sum to $2e$;
and if $j\ge a+3$, the power alone is at least $4e$.
For $w=2^j$ and $t=j$, the bound in \eqref{eq:weightbound},
after subtracting $C_n$ and dividing by $e$, is
\[
 2^{j-a}+a-j+\max\{0,a+2-j\}-1.
\]
It equals $1$ at $j=a+1,a+2$ and is at least $2$ at every other
integer $j\ge0$. Decreasing $t$ adds at least $e$ to the bound.
Thus every term of valuation below $C_n+3e/2$ has
\[
 (w,t)=(2e,a+1)\quad\hbox{or}\quad(w,t)=(4e,a+2).
\]
Both force a single nonzero component $\alpha_i=w$ with $c_i=1$.
If $i\le2^{b-1}$, the actual bound using $I$ raises both contributions
to at least $C_n+2e$. Thus only $2^{b-1}<i\le p^b$ with $c_i=1$ remain.
For each such $i$, the two contributions in \eqref{eq:remcoefficient} are
\begin{equation}\label{eq:binarypair}
 \frac{N}{2e}\left(\binom{N-1}{2ei-1}+1\right)\pi^{2e}x_i^{2e}
       +\frac{N}{4e}\pi^{4e}x_i^{4e}.
\end{equation}
The first binomial coefficient is $1\pmod4$. In its product expression,
the unique factor not congruent to $1\pmod4$ before the overall sign is
$j=N/2$; it is $-1$, and $2ei-1$ is odd. Hence the expression in parentheses,
divided by $2$, is odd.

Modulo $\pi^{e/2}$, both $x_i^{2e}$ and $x_i^{4e}$ equal the representative
$\overline{x_i}\in\{0,1\}$. For a unit $x_i=1+\pi z$, the binomial estimate
$\vp(x_i^{2e}-1)\ge2e$ proves this; for a nonunit it is immediate.
After division by $\pi^{C_n+e}$, \eqref{eq:binarypair} is therefore
\[
 u^{n-a-2}(u^2+1)\overline{x_i}\pmod{\pi^{e/2}}.
\]
The errors from reducing the integer coefficients have valuation at
least $e$, so they vanish modulo $\pi^{e/2}$.
Summing these pairs and the omitted terms proves \eqref{eq:binaryerror}.
\end{proof}

We next evaluate the minimum contributed by the logarithmic term.
For $p>2$, recall the data defining \eqref{eq:rho}: $c=-\eta^{-1}$,
$\ell=\ord(c)$, $L/K$ is unramified of degree $\ell$, and $\zeta$
is the Teichm\"uller lift of a root of $T^{p-1}-c$.
The residue of $\beta\in\OO_L$ generates a normal basis.
For $p=2$, put $c=1$, $\ell=1$, $L=K$, and $\zeta=\beta=1$. We still have
$\delta_0=p^a-ea$ and $C_n=en+\delta_0$.
Let $\sigma$ generate $\operatorname{Gal}(L/K)$, inducing $x\mapsto x^p$ on the residue
field. If $\widetilde c$ is the Teichm\"uller lift of $c$, then
\begin{equation}\label{eq:eigen}
 \sigma(\zeta)=\widetilde c\zeta,\qquad
 \zeta^{p^j}=\widetilde c^{\,j}\zeta.
\end{equation}
Define $\delta$ by
\begin{equation}\label{eq:logminimum}
 \delta=\min_{\varepsilon\in1+\pi\OO_L}
       \vp(\Tr_{L/K}(\zeta\log \varepsilon))\quad(p>2),
 \qquad
 \delta=\min_{\varepsilon\in1+\pi\OO}\vp(\log \varepsilon)\quad(p=2).
\end{equation}
The next lemma evaluates $\delta$ and relates it to $\rho_K$.

\begin{lemma}\label{lem:units}
For odd $p$,
\begin{equation}\label{eq:odddelta}
 \delta=\min\{\delta_0+p^a,\ \vp(\Tr_{L/K}(\zeta\log(1+\pi\beta)))\}.
\end{equation}
For $p=2$, $a\ge1$,
\begin{equation}\label{eq:binarydelta}
 \delta=\min\{\delta_0+3e/2,\ \vp(\log(1+\pi))\}.
\end{equation}
For $p=2$, $a=0$, one has $K=\QQ_2$ and $\delta=2$.
In every case $\delta=\delta_0+\rho_K$, with $\rho_K$ given by
\eqref{eq:rho}, and $\rho_K>0$.
\end{lemma}
\begin{proof}
Suppose $p>2$. The logarithmic terms on $1+\pi^2\OO_L$ have valuations
at least $2j-e v_p(j)$, whose unique minimum is
$2p^a-ea=\delta_0+p^a$ at $j=p^a$. The leading field trace is a nonzero
multiple of
$\Tr_{\FF_{p^\ell}/\FF_p}(\overline\zeta\,\overline x^{p^a})$.
The finite-field trace pairing and Frobenius bijectivity show that it is
nonzero for some $x$. Thus the minimum on $1+\pi^2\OO_L$ is exactly
$\delta_0+p^a$.

Put $\varepsilon_0=1+\pi\beta$. The normal basis condition says that
$\sigma^j(\varepsilon_0)$, $0\le j<\ell$, generate
$(1+\pi\OO_L)/(1+\pi^2\OO_L)$ as an abelian group. By \eqref{eq:eigen},
\[
 \Tr_{L/K}(\zeta\log\sigma^j(\varepsilon_0))
       =\widetilde c^{-j}\Tr_{L/K}(\zeta\log \varepsilon_0).
\]
Every principal unit is a product of integral powers of these generators
and an element of $1+\pi^2\OO_L$. The logarithm is additive on this group.
The preceding formula and the minimum on $1+\pi^2\OO_L$ prove \eqref{eq:odddelta}.

For $p=2$, $e\ge2$, the image of $1+\pi$ generates
$(1+\pi\OO)/(1+\pi^3\OO)$, a group of order four: its square has depth
exactly two, and its fourth power lies in $1+\pi^3\OO$.
On $1+\pi^3\OO$, the minimum of $3j-e v_2(j)$ is uniquely attained at
$j=2^{a-1}$ and equals $\delta_0+3e/2$. Its leading coefficient is a unit
when the parameter is a unit. Additivity gives \eqref{eq:binarydelta}.
For $K=\QQ_2$, principal units are $\{\pm1\}(1+4\ZZ_2)$, the logarithm
kills $-1$, and it is an isometry on $1+4\ZZ_2$; hence $\delta=2$.

For odd $p$, apply Lemma~\ref{lem:error} in dimension one at
$n=a+1$, $b=0$. It gives
\[
 \Tr_{L/K}\bigl(\zeta(\varepsilon_0^{p^{a+1}}-1)\bigr)
 \equiv p^{a+1}\Tr_{L/K}(\zeta\log \varepsilon_0)
              \pmod{\pi^{p^{a+1}+e}}.
\]
Since $p^a<e$, taking the minimum with $p^a$ after subtracting $p^{a+1}$ gives \eqref{eq:rho}.

For $p=2$, $a\ge1$, put $\varepsilon_0=1+\pi$. The binomial expansion gives
$\vp(\varepsilon_0^{2^{a+1}}-1)\ge2^{a+1}>e$ (in fact the two lowest residues
cancel). The logarithm preserves valuation at depths greater than $e$:
if $v_\pi(z)>e$, then for every $j\ge2$ one has
$jv_\pi(z)-e v_2(j)>v_\pi(z)$, so the first term of $\log(1+z)$
has strictly smaller valuation than all later terms.
Therefore
\[
 \vp(\varepsilon_0^{2^{a+1}}-1)=e(a+1)+\vp(\log \varepsilon_0),
\]
including the case where both sides are infinite. This proves
\eqref{eq:rho} for $a\ge1$. For $p=2,a=0$, $\pi$ is a uniformizer
of $\QQ_2$, and $\vp((1+\pi)^2-1)\ge3$, so \eqref{eq:rho}
gives $\rho_K=1$, consistently with $\delta=2$.

The positivity of $\rho_K$ follows from cancellation between the
terms of valuation $\delta_0$. The only terms of this valuation have
indices $p^a,p^{a+1}$. For odd $p$, after normalization their weighted
trace on the residue field is a nonzero scalar multiple of
\[
 \eta\Tr(\overline\zeta x^{p^a})
                +\Tr(\overline\zeta x^{p^{a+1}})=0.
\]
Indeed moving Frobenius through the trace replaces $\overline\zeta$ by
$c^{-1}\overline\zeta$, and $c^{-1}=-\eta$.
For $p=2$, the two residues are equal and add to zero in $\FF_2$.
All other terms have strictly higher valuation. Thus $\delta\ge\delta_0+1$.
\end{proof}

The description of $\delta$ in \eqref{eq:logminimum} also shows that
\eqref{eq:rho} is independent of the auxiliary choices. Normal basis generators give the same minimum by the
argument above. The ratio of any two roots of $T^{p-1}-c$ lies in $\FF_p^\times$.
Its Teichm\"uller lift is a unit in $\QQ_p$, so multiplying
$\zeta$ by this unit does not change the minimum.
If the uniformizer is replaced by $v\pi$, then
$\overline{p/(v\pi)^e}=\eta\overline v^{-e}=\eta$, because
$e=(p-1)p^a$ and the residue field is $\FF_p$.
The principal-unit group itself is unchanged. Thus $\delta$, and hence
$\rho_K=\delta-\delta_0$, is intrinsic to $K$.

For $p=2$ and $a\ge1$, the preceding calculation also controls the
remainder in Lemma~\ref{lem:error}.

\begin{lemma}\label{lem:binarycontrol}
Suppose $p=2$, $f=1$, $a\ge1$, and $M$ satisfies the assumptions of
Lemma~\ref{lem:error}. Then
\[
 \vp(\Tr(M^{2^n})-d)\ge en+\delta=C_n+\rho_K.
\]
\end{lemma}
\begin{proof}
Since $\det M\in1+\pi\OO$, it suffices to show that the remainder
in Lemma~\ref{lem:error} has valuation at least $en+\delta$.
If $\rho_K\le e/2$, its basic bound $C_n+e$ suffices.
If $\rho_K>e/2$, \eqref{eq:binarydelta} implies
$\vp(\log(1+\pi))>\delta_0+e/2$.
The two lowest logarithmic terms have indices $e,2e$; all others have
valuation at least $\delta_0+e/2$. Consequently
\[
 \pi^{-\delta_0}\log(1+\pi)
          \equiv-u^{-a-1}(u+1)\pmod{\pi^{e/2}},
\]
so $\vp(u+1)\ge e/2$. As
$u^2+1=(u+1)^2-2u$, it follows that $\vp(u^2+1)\ge e$.
Lemma~\ref{lem:binaryerror} now gives
$\vp(R_N(M))\ge C_n+3e/2\ge C_n+\rho_K$.
\end{proof}

\Needspace{9\baselineskip}
\begin{proof}[Proof of Theorem~\ref{thm:nonexact}]
Assume the exceptional condition in Theorem~\ref{thm:nonmain}, and retain
$e=(p-1)p^a$, $n=a+b+1$, and $\delta_0=p^a-ea$, so that
$C_n=en+\delta_0$.
Write $\barA\sim\barA_0\oplus\barA_\times$, where $\barA_0$ is
its nilpotent primary part and $\barA_\times$ is its complementary
invertible part.
Proposition~\ref{prop:primary} gives the minimum over the distinct
primary parts. We first compute the common value for every nonzero part.

Let $g\ne T$ be an irreducible primary factor. Its roots satisfy
$T^{p-1}-c=0$.
Each root has degree $\ell=\ord(c)$, since
$\alpha^{p^j}=c^j\alpha$. The ratio of any two roots lies in
$\FF_p^\times$, so their Teichm\"uller lifts have ratio in
$\ZZ_p^\times$. Thus $\delta$ is the same for every nonzero primary factor.
Let $m_g$ be the dimension of this primary part and put $d=m_g/\ell$.
For an arbitrary integral primary lift $A$, work in the finite,
complete commutative algebra
$R=\OO[A]\subseteq\operatorname{End}_{\OO}(\OO^{m_g})$.
Modulo $\pi$, this algebra has a nilpotent ideal with residue field
$\FF_{p^\ell}$. Hensel lifting in $R$ gives the unique lift $Z$ of the
semisimple part of $\barA$ satisfying $Z^{p^\ell}=Z$; equivalently,
$Z=\lim_{j\to\infty}A^{p^{\ell j}}$. The derivative of
$T^{p^\ell}-T$ is $-1$ modulo the maximal ideal, so the lift is unique.
The algebra $\OO[Z]$ is the unramified integer ring $\OO_L$, with
$Z$ corresponding to $\zeta$. Its action commutes with $A$.
The lattice $\OO^{m_g}$ is a finite torsion-free, hence free, $\OO_L$-module of rank $d$.
With respect to an $\OO_L$-basis, write $A=\zeta M$,
where $\overline M$ is unipotent and has the prescribed Jordan block
sizes, all at most $p^b$.
Multiplication by the nonzero scalar $\alpha^{-1}$ does not change the
lengths of the nilpotent Jordan blocks. Finally, the trace over $K$ is
$\Tr_{L/K}$ of the matrix trace over $L$.

If $p>2$, Lemma~\ref{lem:error} and $\rho_K\le p^a<e$ imply
\[
 \Tr_{L/K}\bigl(\zeta^{p^n}(\Tr(M^{p^n})-d)\bigr)
 \equiv p^n\Tr_{L/K}(\zeta^{p^n}\log\det M)
                         \pmod{\pi^{C_n+e}}.
\]
The right side has valuation at least $en+\delta=C_n+\rho_K$.
For $p=2$, all nonzero residual eigenvalues are $1$, so $L=K$;
Lemma~\ref{lem:binarycontrol} gives the same bound when $a\ge1$.
For $p=2,a=0$, Lemma~\ref{lem:error} gives a remainder of valuation
at least $C_n+1=n+2=en+\delta$, since $\delta=2$. Thus it also gives
the required lower bound in this case.
At $n+1$ the same arguments apply with $b+1$. Their bound increases by $e$.
The constant terms $d\Tr_{L/K}(\zeta^{p^n})$ agree at consecutive powers.
It follows that every lift of the nonzero primary type satisfies
$\vp(D_r(A))\ge C_n+\rho_K$.

To attain equality, choose $\varepsilon\in1+\pi\OO_L$ attaining $\delta$ in
Lemma~\ref{lem:units}. If $\delta$ is attained at $1+\pi\beta$ for odd $p$,
or at $1+\pi$ for $p=2$, use that unit. Otherwise use $\varepsilon=1+\pi^2x$ for odd $p$, choosing its leading trace
nonzero, or $\varepsilon=1+\pi^3$ for $p=2,a\ge1$. For $p=2,a=0$, use
$\varepsilon=1+\pi^2$, whose logarithm has valuation $2$. Change only one diagonal entry
in one upper triangular block from $\zeta$ to $\zeta \varepsilon$.
The change in the trace at exponent $p^n$ is
\[
 \Tr_{L/K}\bigl(\zeta^{p^n}(\varepsilon^{p^n}-1)\bigr).
\]
There is no factor equal to the block length: only one diagonal entry
has been changed. Also $\zeta^{p^n}=\widetilde c^{\,n}\zeta$, so the
minimum in \eqref{eq:logminimum} is unchanged.
For odd $p$, the scalar form of Lemma~\ref{lem:error} shows its valuation
is exactly $en+\delta$, because the error is strictly higher.
For $p=2$, $\varepsilon^{2^n}$ has depth greater than $e$; the isometry of the
logarithm gives valuation exactly $en+\delta$ for $\varepsilon^{2^n}-1$.
The change at exponent $p^{n+1}$ has strictly larger valuation.
For $p=2$, the valuation of the trace difference also follows directly from
$\varepsilon^{2^{n+1}}-\varepsilon^{2^n}=\varepsilon^{2^n}(\varepsilon^{2^n}-1)$.
The unchanged triangular blocks have zero trace difference. We have proved
\begin{equation}\label{eq:nonzeroexact}
 E_r(K,\barA_\times)=C_n+\rho_K.
\end{equation}
Together with \eqref{eq:primarymin}, this proves the exact formula with
$E_r(K,\barA_0)$ in place of $L_{h_0,e}(n)$.

Finally, we show that $E_r(K,\barA_0)$ can be replaced by $L_{h_0,e}(n)$
in the minimum. By Theorem~\ref{thm:complete}, the only exception for odd $p$
requires $e=(p-1)p^v/2$, incompatible with $e=(p-1)p^a$.
For $p=2$, its case~\textup{(II)} requires $e=3\cdot2^v$, also incompatible.
Only its case~\textup{(I)} can remain: $e=2^a$, $h_0=3\cdot2^s$, with
$\nu=n-s-a\ge6$ even and further partition restrictions.
Under these conditions, \eqref{eq:Lvalues} gives
\[
 L_{h_0,e}(n)=es+e(2^\nu/3+\nu-1/3),\qquad
 C_n=e(s+\nu+1).
\]
Hence $L_{h_0,e}(n)-C_n=e(2^\nu-4)/3\ge20e>\rho_K$.
Whether or not its nilpotent correction is present, that component cannot
be the minimum in \eqref{eq:primarymin}. In all other nilpotent cases its
exact value is $L_{h_0,e}(n)$. This proves \eqref{eq:nonanswer} and attainment.
The strict inequality follows from Corollary~\ref{cor:improvement};
the positivity and independence of $\rho_K$ were established in
Lemma~\ref{lem:units} and the discussion following its proof.
\end{proof}

\section{Dependence on the local field}\label{sec:arithmetic}
Outside the three nilpotent exceptions, $E_r$ is determined by
$p,e,n,\lambda_1$. For nonnilpotent classes not satisfying the
exceptional condition in Theorem~\ref{thm:nonmain}, it equals $C_n$
and is determined by $p,e,n$. In the exceptional cases, the formulas
also involve the unit $u=p/\pi^e$ or principal units of $K$.

Recall that the nonnilpotent exceptional case has $f=1$ and
$e=(p-1)p^a$. With $L$ and $\zeta$ as in Lemma~\ref{lem:units},
\eqref{eq:logminimum} and Lemma~\ref{lem:units} give
\begin{equation}\label{eq:arithmeticlog}
 C_n+\rho_K=en+
 \min_{\varepsilon\in1+\pi\OO_L}\vp\!\left(\Tr_{L/K}(\zeta\log \varepsilon)\right).
\end{equation}
The evaluation in Lemma~\ref{lem:units} yields
\begin{equation}\label{eq:rhorange}
 1\le\rho_K\le
 \begin{cases}p^a,&p>2,\\1,&p=2,\ a=0,\\3\cdot2^{a-1},&p=2,\ a\ge1.
 \end{cases}
\end{equation}
These upper endpoints occur over cyclotomic fields. Let
$K=\QQ_p(\xi)$, where $\xi$ is a primitive $p^{a+1}$st root of unity,
and take $\pi=\xi-1$. The polynomial $\Phi_{p^{a+1}}(1+T)$ is
Eisenstein, of degree $(p-1)p^a$, so the field has the required
ramification index and residue field. Its constant term is $p$ and
all intermediate coefficients are divisible by $p$. Substitution of
$\pi$ and division by $\pi^e$ therefore give $\eta=-1$.
We may take $L=K$, $\zeta=\beta=1$ in \eqref{eq:rho}; since
$(1+\pi)^{p^{a+1}}=1$, $\rho_K$ attains the corresponding upper bound
in \eqref{eq:rhorange}.
For the scalar residual class $(1)$ and $n\ge a+1$, the exponent is
$C_n+\rho_K$. Thus the endpoints in \eqref{eq:rhorange} cannot be decreased uniformly
over fields with the indicated parameters.

When $p>2$ and $\eta=-1$, we have $c=1$ and may take $L=K$,
so \eqref{eq:rho} becomes
\[
 \rho_K=\min\{p^a,\ \vp((1+\pi)^{p^{a+1}}-1)-p^{a+1}\}.
\]

The following examples show that fields with the same $(p,e,f)$
can give different exponents for the same residual class.

\begin{example}\label{ex:binaryfields}
Let $K_1=\QQ_2(i)$ with $\pi_1=i-1$, and
$K_2=\QQ_2(\sqrt2)$ with $\pi_2=\sqrt2$.
Both have $(p,e,f)=(2,2,1)$, so $a=1$.
Since $(1+\pi_1)^4-1=0$, formula \eqref{eq:rho} gives $\rho_{K_1}=3$.
On the other hand,
$(1+\pi_2)^4-1=16+12\pi_2$ has valuation $5$:
its two terms have valuations $8$ and $5$, respectively.
Thus $\rho_{K_2}=1$. For the same scalar residual class $(1)$,
\[
 E_{n+1}(K_1,(1))=2n+3,\qquad
 E_{n+1}(K_2,(1))=2n+1\qquad(n\ge2).
\]
\end{example}

\begin{example}\label{rem:field}
Take $p=q=2$, $e=4$, $n=8$, and residual type $J_3(0)$. Then
the minimum \eqref{eq:introL} is $L_{3,4}(8)=108$.
In $K_1=\QQ_2(\pi_1)$ with $\pi_1^4=2$, the unit $u=2/\pi_1^4$ is $1$,
and \eqref{eq:answer} gives $E_9=112$.
In $K_2=\QQ_2(\pi_2)$, where $\pi_2$ is a root of the Eisenstein
polynomial $T^4+2T^3+2T^2+2$, put $B=1+\pi_2^2+\pi_2^3$.
Then $\pi_2^4=-2B$ and $u=2/\pi_2^4=-B^{-1}$. Modulo $\pi_2^8$,
\[
 B^4(u^2+1+\pi_2^6u^4)=B^2+B^4+\pi_2^6
 \equiv4\pi_2^2+4\pi_2^3+2\pi_2^5+2\pi_2^6\equiv0.
\]
Thus \eqref{eq:answer} gives $E_9=116$, whereas $E_9=112$ over $K_1$.
Both fields have $(p,e,f)=(2,4,1)$.
\end{example}

Examples~\ref{ex:binaryfields} and \ref{rem:field} show that, even
for fixed $r$ and residual similarity class, $E_r$ is not determined
by $(p,e,f)$. Although the formulas use a uniformizer, their values
are independent of its choice: by Theorems~\ref{thm:complete} and
\ref{thm:nonexact}, they equal the minimum in \eqref{eq:definition}.

\end{document}